\documentclass[a4paper,10pt,leqno]{article}
\usepackage{amsfonts,amssymb,amsmath,amsthm}
\usepackage{fullpage}
\usepackage{epic,eepic}
\usepackage{mathrsfs}
\newtheorem{theorem}{Theorem}[section]
\newtheorem{proposition}[theorem]{Proposition}

\newtheorem{lemma}[theorem]{Lemma}

\theoremstyle{remark}
\newtheorem{remark}[theorem]{{\bf Remark}}

\numberwithin{equation}{section}

\newcommand{\A}{\mathbf A}
\newcommand{\N}{{\mathbb N}}
\newcommand{\R}{{\mathbb R}}

\newcommand{\eps}{{\varepsilon}}

\renewcommand\a{\alpha}
\renewcommand\b{\beta}
\newcommand\g{\gamma}

\renewcommand\d{\delta}
\renewcommand\l{\lambda}
\newcommand\n{\nabla}
\renewcommand\o{\omega}\renewcommand\O{\Omega}

\def\var{\varphi}

\def\vark{\varkappa}

\def\calF{{\cal F}}

\def\calF{{\cal F}}

\newcommand{\esssup}{\mathop{\mathrm{ess\,sup}}}
\newcommand{\essinf}{\mathop{\mathrm{ess\,inf}}}
\newcommand{\essosc}{\mathop{\mathrm{ess\,osc}}}
\renewcommand{\div}{\operatorname{div}}

\newcommand{\ind}{\hbox{\rm 1\hskip -4.5pt 1}}

\def\part{\partial}

\def\pr{^{\prime}}

\def\XXint#1#2#3{{\setbox0=\hbox{$#1{#2#3}{\int}$}
\vcenter{\hbox{$#2#3$}}\kern-.5\wd0}}

\def\W{\rlap{$\buildrel \circ \over W$}\phantom{W}}
\def\V{\rlap{$\buildrel \circ \over V$}\phantom{V}}

\begin{document}

\title{
Harnack inequality and continuity
of solutions
 to quasi-linear degenerate parabolic equations with coefficients from Kato-type classes }
\author{
{\Large Vitali Liskevich}\footnote{Corresponding author}\\
\small Department of Mathematics\\
\small University of Swansea\\
\small Swansea SA2 8PP, UK\\
{\tt v.a.liskevich@swansea.ac.uk}\\
\and
{\Large Igor I.\,Skrypnik}\\
\small Institute of Applied\\
\small Mathematics and Mechanics\\
\small Donetsk 83114,  Ukraine\\
{\tt iskrypnik@iamm.donbass.com} }

\date{}

\maketitle

\setlength{\unitlength}{0.0004in}
\begingroup\makeatletter\ifx\SetFigFont\undefined%
\gdef\SetFigFont#1#2#3#4#5{%
  \reset@font\fontsize{#1}{#2pt}%
  \fontfamily{#3}\fontseries{#4}\fontshape{#5}%
  \selectfont}%
\fi\endgroup%
\renewcommand{\dashlinestretch}{30}

\begin{abstract}
For a general class of divergence type quasi-linear degenerate parabolic  equations with
measurable coefficients and lower order
terms from non-linear Kato-type classes,
we prove local boundedness and continuity of solutions, and the intrinsic
Harnack inequality for positive solutions.
\end{abstract}

\bigskip
Keywords: quasi-linear degenerate parabolic  equations, local boundedness, continuity, Harnack inequality.

\section{Introduction and main results}

In this paper we are concerned with  general divergence type quasi-linear degenerate parabolic  equations with
measurable coefficients and lower order
terms. This class of equations has numerous applications and has been attracting attention for several decades (see, e.g.  the monographs
\cite{DiB, LaSU, WUZ}, survey \cite{DiUV} and references therein).

Let $\O$ be a domain in $\R^n$, $T>0$. Set $\O_T=\O\times (0,T)$. We study
solutions to
the
equation
\begin{equation}
\label{e1.1} u_t-\div \A (x,t,u, \n u)=a_0(x,t,u,\n u) ,\quad (x,t)\in
\O_T.
\end{equation}

Throughout the paper we suppose that 
the functions $\A:\O\times\R^+\times\R\times \R^n\to \R^n$ and
$a_0:\O\times \R^+\times\R\times \R^n\to \R^n$ are such that
$\A(\cdot,\cdot,u,\zeta),\ a_0(\cdot,\cdot,u,\zeta)$ are Lebesgue measurable for all $u\in \R,\zeta\in\R^N$,
and $\A(x,t,\cdot,\cdot),\ a_0(x,t,\cdot,\cdot)$ are continuous for almost all $(x,t)\in\O_T$.

We also assume that the following structure conditions are satisfied:
\begin{eqnarray}
\nonumber
\A (x,t,u, \zeta)\zeta &\ge& c_1 |\zeta|^p,\quad \zeta\in\R^n,\\
|\A (x,t,u, \zeta)|&\le& c_2|\zeta|^{p-1}+g_1(x) |u|^{p-1}+f_1(x),\label{e1.2b}\\
\nonumber |a_0(x,t,u,\zeta)|&\le& h(x)|\zeta|^{p-1}+g_2(x)|u|^{p-1}+f_2(x),
\end{eqnarray}
where $2<p<n$,
$c_1, c_2$ are positive constants and $f_1(x),f_2(x),g_1(x),g_2(x), h(x)$ are nonnegative functions,
satisfying conditions which will be specified below. The constants in \eqref{e1.2b}, $n$ and $p$ are further referred to as the data.
The aim of this paper is to establish basic qualitative properties such as local boundedness
of weak solutions, their continuity and the Harnack inequality for positive solutions under minimal
possible restrictions on the coefficients in  structure conditions~\eqref{e1.2b}.
These properties are indispensable in the qualitative
theory of second-order elliptic and parabolic equations. For equation \eqref{e1.1} with $g_1=g_2=h=0$ and $f_1,f_2$ constants
the local boundedness and H\"older continuity of solutions was know  since mid-1980s (see \cite{DiB,DiUV} for the results, references and historical notes),
and a recent break through has been made in \cite{DiGV}, where the intrinsic Harnack inequality has been proved.
Before stating precisely our results we
make several remarks related to lower order terms of \eqref{e1.1} and
refer the reader  for an extensive survey of the regularity issues
to \cite{DiB,DiUV,DiGV}.



Local boundedness and H\"older continuity of weak solutions to
homogeneous linear divergence type second-order elliptic equations with measurable coefficients
without lower order terms is known since the famous results by De Giorgi \cite{degiorgi} and Nash \cite{N},
and the Harnack inequality since
 Moser's  celebrated paper \cite{Mo}. However in presence
of lower order term in the equation weak solutions may have singularities and/or internal zeroes, and
the Harnack inequality in general may not be valid, as one can easily realise looking
at the equation $-\Delta u +\frac{c}{|x|^2}u=0$. It was Serrin~\cite{serrin} who generalized Moser's result to the case of
 quasi-linear equations with lower order
terms with conditions expressed in terms of  $L^p$-spaces. 
Using probabilistic techniques Aizenman and Simon in their famous paper~\cite{AS} proved the Harnack inequality and continuity of weak solutions
to the equation $-\Delta u +V u=0$ under the local Kato class condition on the potential $V$. Moreover, they showed that the Kato type condition on the potential $V$ is necessary
for the validity of the Harnack inequality.
Soon after that Chiarenza, Fabes and Garofalo~\cite{CFG} developed
a real variables techniques to prove the Harnack inequality for a linear equation of divergence
type with measurable coefficients and the potential from the Kato class,
thus extending Aizenman, Simon's result. Kurata~\cite{kurata} extended the method of Chiarenza, Fabes and Garofalo and
proved the same for the equation $-\sum_{k,j}\partial_k a_{kj}\partial_j u+\sum_{j}b_j\partial_j u+Vu=0$.
with $|b|^2,V$ from the Kato class. Both papers \cite{CFG} and \cite{kurata} make a heavy use of Green's functions which makes this approach inapplicable
to quasi-linear equations.
To treat the quasi-linear case of $p$-Laplacian with a lower order term Biroli \cite{biroli0,biroli} introduced the notion
of the nonlinear Kato class and gave the Harnack inequality
for positive solutions to $-\Delta_p u +V u^{p-1}=0$.
This was extended in \cite{Skr2} to the general case of quasi-linear elliptic equations with
lower order terms.

For second-order linear parabolic equations with measurable coefficients  (without lower order terms) H\"older continuity of solutions was first proved by Nash \cite{N}.
Moser~\cite{Mo1} proved the validity of the Harnack inequality which was extended to the case of quasilinear equations with $p=2$
in the structure conditions and structure coefficients from $L^q$-classes in \cite{ArSe}, \cite{T}.
The continuity of weak solutions and the Harnack inequality for second-order linear elliptic equations with lower order coefficients
from Kato-classes was proved by Zhang \cite{Z,Z1}.

The parabolic theory for degenerate quasi-linear equations differs substantially from the "linear" case $p=2$ which can be already realized looking at the Barenblatt solution
to the parabolic $p$-Laplace equation.
DiBenedetto developed an innovative intrinsic scaling method (see \cite{DiB} and the references to the original papers there; see also a nice exposition in \cite{urbano} where some recent advances are included)
and proved the H\"older continuity of weak solutions to \eqref{e1.1}
for $p\not=2$ for the case $g_1=g_2=h=0$ and $f_1,f_2$ from $L^q$-classes, and the intrinsic Harnack inequality for the parabolic $p$-Laplace equations. For the measurable coefficients in the main
part of \eqref{e1.1} the intrinsic Harnack inequality was proved in the recent break-through paper \cite{DiGV}. It is natural to conjecture that the intrinsic Harnack inequality holds
for the parabolic $p$-Laplace equation perturbed by lower order terms with coefficients from Kato classes. The difficulty is that seemingly neither De Giorgi nor Moser iteration
techniques work in this situation.



In this paper following the strategy of \cite{DiGV} but using a different iteration, namely  the Kilpel\"ainen-Mal\'y technique \cite{KiMa}
properly adapted to the parabolic equations (cf.\,\cite{Skr,Skr1}),
we establish the local boundedness and continuity of solutions to \eqref{e1.1} and the intrinsic Harnack inequality.

Following Biroli \cite{biroli0,biroli} we introduce the non-linear Kato $K_p$ class by
\begin{equation}
\label{e1.9}
K_p:=\left\{g\in L^1(\O)\,:\,
\lim_{R\to 0}\sup_{x\in\O}
\int_0^R\left\{\frac{1}{r^{n-p}}
\int_{B_r(x)\cap \O}|g(z)|dz
\right\}^{\frac{1}{p-1}}\frac{dr}{r}=0
\right\},
\end{equation}
where $B_r(x)=\{z\in \O\,:\,|z-x|<r\}$.
As one can easily see, for $p=2$, $K_p$ reduces to the standard definition of the Kato class as defined in \cite{AS,simon}.

We will also need the class $\widetilde{K_p}$ of functions $g\in
L^1(\O)$ satisfying the condition
\begin{equation}
\label{e2.6} \lim_{R\to 0}
\sup_{x\in\O}\int_0^R\left\{\frac{1}{r^{n-p}}
\int_{B_r(x)}|g(z)|dz\right\}^{\frac{1}{p}}\frac{dr}{r}=0.
\end{equation}
It is easy to see that $\widetilde{K_p}\subset K_p$.
We assume that
\begin{eqnarray}
\label{e1.6b}
 F_1:=(g_1+f_1)^\frac{p}{p-1}\in \widetilde{K_p},\\
 \label{e1.7b}
F_2:=h^p+g_2+f_2\in K_p.
\end{eqnarray}

In what follows we use the following quantities
\begin{eqnarray}
\label{e1.8b}
\calF_1(R)=\sup_{x\in\O}\int_0^R\left(\frac1{r^{n-p}}\int_{B_r(x)}F_1(z)dz\right)^\frac1{p}\frac{dr}{r},\\
\label{e1.9b}
\calF_2(R)=\sup_{x\in\O}\int_0^R\left(\frac1{r^{n-p}}\int_{B_r(x)}F_2(z)dz\right)^\frac1{p-1}\frac{dr}{r}.
\end{eqnarray}

\bigskip

Before formulating the main results, let us remind the reader of the definition of a weak solution
to equation \eqref{e1.1}.

We say that $u$ is a weak solution to \eqref{e1.1} if $u\in
V(\O_T):=W^{1,p}_{loc}(\O_T)\cap C([0,T]; L^2_{loc}(\O))$ and for
any interval $[t_1,t_2]\subset (0,T)$ the integral identity
\begin{equation}
\label{e1.12b}
\int_{\O}u\var dx\Big|_{t_1}^{t_2}+\int_{t_1}^{t_2}
\int_\O \left\{-u\var_\tau+\A(x,\tau,u,\n u) \n\var -a_0(x,\tau,u,\n u) \var\
\right\}dx\,d\tau=0
\end{equation}
for any $\var \in \V(\O_T)$.

\medskip

The first main result of this paper  is the local boundedness of
solutions.

\begin{theorem}
\label{thm1.1b}
 Let conditions \eqref{e1.2b}, \eqref{e1.6b} and \eqref{e1.7b}
be fulfilled. Let $u$ be a weak solution to equation~\eqref{e1.1}.
Then $u$ locally bounded, that is $u\in L^\infty_{loc}(\O_T)$.
\end{theorem}

The proof of Theorem~\ref{thm1.1b} is based on the adaptation of the
Kilpel\"ainen-Mal\'y technique \cite{KiMa} to parabolic equations
using ideas from \cite{Skr,Skr1}.
Having established the local boundedness we proceed with the
continuity. At this stage we can assume that the solutions are
bounded in $\O_T$.

\begin{theorem}
\label{thm1.1}
Let conditions \eqref{e1.2b}, \eqref{e1.6b} and
\eqref{e1.7b} be fulfilled and $h(x)=1$. Let $u$ be a bounded weak solution to
equation~\eqref{e1.1}. Then $u$ is continuous, that is $u\in C(\O_T)$.
\end{theorem}

Next is the Harnack inequality for positive solutions to \eqref{e1.1}.

Let $u$ be a nonnegative solution to \eqref{e1.1}.
Fix a point $(x_0,t_0)\in \O_T$ such that $u(x_0,t_0)>0$. Consider the cylinders
\[
Q_\rho^\theta(x_0,t_0)=B_\rho(x_0)\times (t_0-\theta \rho^p,t_0+\theta \rho^p), \quad \theta=\left(\frac{c}{u(x_0,t_0)}\right)^{p-2},
\]
where $c>0$ is fixed.
\begin{theorem}
\label{Harnack}
Let conditions the conditions of Theorem~\ref{thm1.1}
 be fulfilled and $h(x)=1$.
Let $u$ be a positive solution to \eqref{e1.1}.
Then there exist positive constants $c,\g$ depending only on the data and $\sup_{\O_T}u(x,t)$, such that
for all intrinsic cylinders $Q_{4\rho}^\theta(x_0,t_0)\subset \O_T$
either $u(x_0,t_0)\le \g (\rho +   \calF_1(2\rho)+\calF_2(2\rho))  $
or
\[
u(x_0,t_0)\le \g \inf_{B_\rho(x_0)} u(x,t_0+\theta \rho^p),\quad \theta=\left(\frac{c}{u(x_0,t_0)}\right)^{p-2}.
\]
Moreover, if $g_1=g_2=0$, the constant $\g$ can be chosen independent from $\sup_{\O_T}u(x,t)$.
\end{theorem}

\bigskip

\begin{remark}
In the linear theory the Kato class is known to be the optimal condition on the zero order term of the equation $\Delta u +Vu=0$
to imply the continuity of solutions and the Harnack inequality. The same is true for the quasi-linear equations. For the equation
$\Delta_p u + V u |u|^{p-2}=0$ with $V$ behaving around zero like $c\left(\log\frac1{|x|}\right)^{1-p}$, depending on the sign of $V$,
one can easily produce a solution with singularity at zero, or with internal zero at zero (see, e.g.\cite{LS,LS1}).
On the other hand, the function  $c\left(\log\frac1{|x|}\right)^{1-p+\eps}$ is from the Kato class $K_p$ for any $\eps>0$,
and Theorems~\ref{thm1.1b},~\ref{thm1.1},~\ref{Harnack} apply.
\end{remark}

\bigskip

Let $u$ be a weak solution to \eqref{e1.1} in $\O_T$. Let $(y,s)\in \O_T$ be an arbitrary point. Consider
the cylinder $Q_{4\rho}^\theta(y,s)\subset \O_T$,
\[
Q_\rho^\theta(y,s)=B_\rho(y)\times (s-\theta \rho^p,s), \quad \theta>0.
\]

Denote by $\mu_{\pm}$ and $\o$ non-negative numbers such that
\[
\mu_+\ge \esssup_{Q_{4\rho}^\theta(y,s)}u(x,t),\quad \mu_-\le \essinf_{Q_{4\rho}^\theta(y,s)}u(x,t),\quad \o\ge \mu_+-\mu_-.
\]

\bigskip
As was already mentioned, our strategy of the proof of the Harnack inequality is the same as in \cite{DiGV}. Namely,
Theorems~\ref{thm1.1},~\ref{Harnack} will be consequences of the following two theorems.

The next theorem is a De Giorgi-type lemma (cf.\,\cite{DiGV}), and its formulation is almost the same as in   \cite{DiGV}.
However, due to the different structure conditions the De Giorgi type iteration cannot be used. Instead, we adapt the Kilpel\"ainen-Mal\'y iteration \cite{KiMa} combined
with ideas from \cite{Skr,Skr1}, where the Kilpel\"ainen-Mal\'y technique was adapted to parabolic equations.
\begin{theorem}
\label{thm1.4}
Let the conditions of Theorem~\ref{thm1.1} be fulfilled.
Fix $\xi,a\in(0,1)$, $(\xi\o)^{p-2}\ge \frac1{\theta}$.
There exist numbers $B\ge 1$ and $\nu\in(0,1)$ depending only on the data and $\theta,\xi,\o$ and $a$ such that if
\begin{equation}
\label{e1.12}
\left|\{(x,t)\in Q_{2\rho}^\theta(y,s)\,:\,u(x,t)\le \mu_-+\xi\o\}\right|\le \nu | Q_{2\rho}^\theta(y,s)|,
\end{equation}
then either $\xi\o\le B(\rho+\calF_1(2\rho)+\calF_2(2\rho))$, or
\begin{equation}
\label{e1.13}
u(x,t) \ge \mu_-+ a \xi\o \ \ \text{for almost all (a.a.)}\ (x,t)\in  Q_{\rho}^\theta(y,s).
\end{equation}
Likewise,if
\begin{equation}
\label{e1.14}
\left|\{(x,t)\in Q_{2\rho}^\theta(y,s)\,:\,u(x,t)\ge \mu_+-\xi\o\}\right|\le \nu | Q_{2\rho}^\theta(y,s)|,
\end{equation}
then either $\xi\o\le B(\rho+\calF_1(2\rho)+\calF_2(2\rho))$, or
\begin{equation}
\label{e1.15}
u(x,t) \le \mu_+- a \xi\o \ \ \text{for almost all (a.a.)}\ (x,t)\in  Q_{\rho}^\theta(y,s),
\end{equation}
where $\calF_1(\rho),\, \calF_2(\rho)$ are defined in \eqref{e1.6b}, \eqref{e1.7b}.
\end{theorem}

The following theorem is an expansion of positivity result, analogous in formulation as well as in the proof to \cite[Lemma~3.1]{DiGV}.

\begin{theorem}
\label{thm1.5}
Let the conditions of Theorem~\ref{thm1.1} be fulfilled.
There exist positive numbers $B,b_1<b_2$ and $\sigma\in(0,1)$ depending only on the data
such that if
\begin{equation}
\label{e1.16}
u(x,s) \ge \mu_-+N \quad \text{for} \ x\in B_\rho(y),
\end{equation}
then either $N\le B(\rho +\calF_1(2\rho)+\calF_2(2\rho))$, or
\begin{equation}
\label{e1.17}
u(x,t) \ge \mu_-+\sigma N \quad \text{for a.a.} \ x\in B_{2\rho}(y),
\end{equation}
for all
\begin{equation}
\label{e1.18}
s+N^{2-p} b_1 \rho^p\le t\le s+N^{2-p} b_2\rho^p.
\end{equation}
If on the other hand
\begin{equation}
\label{e1.19}
u(x,s) \le \mu_+-N \quad \text{for} \ x\in B_\rho(y),
\end{equation}
then either $N\le B(\rho +\calF_1(2\rho)+\calF_2(2\rho))$, or
\begin{equation}
\label{e1.20}
u(x,t) \le \mu_+-\sigma N \quad \text{for a.a.} \ x\in B_{2\rho}(y),
\end{equation}
for all $t$ satisfying \eqref{e1.18}.
\end{theorem}

The rest of the paper contains the proof of the above theorems. In Section~\ref{aux} we collect some auxiliary propositions and
required integral estimates os solutions. In Section~\ref{bdd} we give a proof of local boundedness of solutions which is based on the parabolic modification of
the Kilpel\"ainen-Mal\'y technique \cite{KiMa}. Section~\ref{degiorgi} contains the proof of the variant of De Giorgi lemma, Theorem~\ref{thm1.4}.
Expansion of positivity, Theorem~\ref{thm1.5}, is proved in Section~\ref{expansion}. In Section~\ref{cont} we prove continuity of solutions following
\cite{DiB}.  Finally, in Section~\ref{harnack} we sketch a proof of the intrinsic Harnack inequality, Theorem~\ref{Harnack}, leaving out details for which we refer
to \cite{DiGV}.

\section{Auxiliary material and integral estimates of solutions}
\label{aux}
\subsection{Local energy estimates}

\begin{lemma}
\label{lem2.2}
Let $u$ be a solution to \eqref{e1.1} in $\O_T$. Then there exists $\g>0$ depending only on $n,p,c_1,c_2$ such that for every cylinder
$Q_\rho^\theta(y,s)=B_\rho(y)\times (s-\theta\rho^p,s)\subset \O_T$, any $k\in\R^1$ and any smooth $\xi(x,t)$
which is zero for $(x,t)\in \partial B_\rho(y)\times (s-\theta\rho^p,s)$ one has
\begin{eqnarray}
\nonumber
&&\sup_{s-\theta\rho^p<t<s}\int_{B_\rho(y)}(u-k)^2_{\pm}\xi(x,t)^pdx + c_1\iint_{Q_\rho^\theta(y,s)}|\n(u-k)_{\pm}|^p\xi(x,t)^pdxdt\\
\nonumber
&\le& \int_{B_\rho(y)}(u-k)^2_{\pm}\xi(x,s-\theta\rho^p)^pdx + \g\iint_{Q_\rho^\theta(y,s)}\left((u-k)^p_{\pm}|\n \xi|^p+(u-k)^2_{\pm}|\xi_t|\right)dxdt\\
\nonumber
&+&\g \iint_{A_\rho(y,s)}[f_1(x)^\frac{p}{p-1}+g_1(x)^\frac{p}{p-1}|u|^p]\xi^pdxdt+
\g \iint_{Q_\rho^\theta(y,s)}(u-k)^p_{\pm}h(x)^p\xi^pdxdt+\\
\label{e2.4}
&+&\g \iint_{Q^\theta_\rho(y,s)}(u-k)_\pm(f_2(x)+g_2(x)|u|^{p-1})\xi^p dxdt,
\end{eqnarray}
where $A_\rho(y,s)=Q_\rho^\theta(y,s)\cap\{(u-k)_{\pm}>0\}$.
\end{lemma}
\proof Test \eqref{e1.1} by $\var=(u-k)_{\pm}\xi^p$ and use conditions \eqref{e1.2b} and the H\"older and Young inequalities.\qed

\bigskip

Let
\[
H^\pm_k:=\esssup_{Q_\rho^\theta(y,s)}|(u-k)_\pm|,\quad \Psi_\pm(u):=\left(\ln\frac{H^\pm_k}{H^\pm_k-(u-k)_\pm+c}\right)_+,\ 0<c<H^\pm_k.
\]
\begin{lemma}
\label{lem2.3}
Let $u$ be a solution to \eqref{e1.1} in $\O_T$. Then there exists $\g>0$ depending only on $n,p,c_1,c_2$ such that for every cylinder
$Q_\rho^\theta(y,s)=B_\rho(y)\times (s-\theta\rho^p,s)\subset \O_T$, any $k\in\R^1$ and any smooth $\xi(x)$
which is zero for $|x-y|>\rho$ one has
\begin{eqnarray}
\nonumber
&&\sup_{s-\theta\rho^p<t<s}\int_{B_\rho(y)}\Psi^2_\pm(u)\xi^pdx
\le \int_{B_\rho(y)\times\{s-\theta\rho^p\}} \Psi^2_\pm(u)\xi^pdx
+\g \iint_{Q_\rho^\theta(y,s)}\Psi_\pm|\Psi_\pm\pr(u)|^{2-p}|\n \xi|^pdxdt\\
\nonumber
&&+ \g \iint_{A_\rho(y,s)}\Psi_\pm|\Psi_\pm\pr(u)|^{2-p}h(x)^pdxdt
+\g\iint_{A_\rho(y,s)}\Psi_\pm|\Psi_\pm\pr(u)|^{2}\left(f_1(x)^\frac{p}{p-1}+g_1(x)^\frac{p}{p-1}|u|^p\right)dxdt\\
\label{e2.5}
&&+\g\iint_{A_\rho(y,s)}\Psi_\pm|\Psi_\pm\pr(u)|\left(f_2(x)+g_2(x)|u|^{p-1}\right)dxdt.
\end{eqnarray}
\end{lemma}
The proof is analogous to that of \cite[Proposition~3.2,\,Chapter~II]{DiB}.

\subsection{Auxiliary propositions}

The following two lemmas will be used in the sequel. The first one is the well known De Giorgi-Poncar\'{e} lemma (see \cite[Chapter~I]{DiB}, \cite[Chapter~II, Lemma~3.9]{La}).

\begin{lemma}
\label{lem2.3m}
Let $u\in W^{1,1}(B_\rho(y))$ for some $\rho>0$ and $y\in \R^n$. Let $k$ and $l$ be real numbers such that $k<l$. Then there exists a constant $\g$
depending only on $n$ such that
\begin{equation}
\label{e2.3m}
(l-k)|A_{k,\rho}||B_\rho\setminus A_{l,\rho}|\le \g \rho^{n+1}\int_{A_{l,\rho}\setminus A_{k,\rho}}| \n u | dx,
\end{equation}
where $A_{k,\rho}=\{x\in B_\rho\,:\,u(x)<k\}$.
\end{lemma}

The next lemma is the time-dependent  version of the measure-theoretic lemma from \cite{DiGV-1}, which can be extracted from
\cite[Section~8]{DiGV}.
\begin{lemma}
\label{lem2.1}
Let $Q_1=B_1(0)\times (-1,0)$ and $v\in V(Q_1)$.
Let $v$ satisfy \eqref{e2.4}.
Suppose that there exist constants $\g>0$ and $\nu\in(0,1)$ such that
\begin{equation}
\label{e2.4m}
\iint_{Q_1} |\n v |^pdxdt\le \g \quad \text{and}\quad |\{(x,t)\in Q_1\,:\,v(x,t)>1\}|>\nu.
\end{equation}
Then for any $\l\in(0,1)$ and $\nu_0\in (0,1)$ there exist a point $(y,s)\in Q_1$, a number $\eta_0\in(0,1)$ and a cylinder $Q_{2\eta_0}(y,s)\subset Q_1$ such that
\begin{equation}
\label{e2.5m}
|\{(x,t)\in Q_{\eta_0}(y,s)\,:\,v(x,t)>\l\}|\ge(1-\nu_0)|Q_{\eta_0}(y,s)|,
\end{equation}
where $Q_R(y,s)=B_R(y)\times (s-R^p,s)$.
\end{lemma}

%
In what follows we will frequently use
the following lemma which is due to Biroli \cite{biroli0,biroli}.

\begin{lemma}
\label{lem2.1b}Let $1<q<n$.
For any $\eps>0$ there exist $R_0<1$ and $\tau>0$ such that the inequality
\begin{equation}
\label{e2.1b}
\sup_{x\in B_1(0)}\int_0^{R_0}\left(\frac1{r^{n-q}}\int_{B_r(x)\cap B_1(0)} H(z)dz\right)^{\frac1{q-1}}\frac{dr}{r}<\tau
\end{equation}
implies that
\begin{equation}
\label{e2.2b}
\int_{B_R(x_0)}H(x)|\var(x)|^qdx \le \eps\int_{B_R(x_0)}|\nabla \var(x)|^qdx
\end{equation}
for any $\var\in \W^{1,q}(B_R(x_0))$ if $R\le R_0$ and $B_{4R_0}(x_0)\subset B_1(0)$.
\end{lemma}

\medskip
\subsection{Integral estimates of solutions}
Set
\[
G(u)=
\left\{
\begin{array}{lll}
u&\text{for}\ &\ u>1,\\
u^{2-2\l}&\text{for}\ &\ 0<u\le 1.
\end{array}
\right.
\]

\begin{lemma}
\label{lem2.2b}
Let the conditions of Theorem~\ref{thm1.1b} be fulfilled. Let $u$ be a solution to \eqref{e1.1}.
Then there exists a constant $\g>0$ depending only on $n,p,c_1,c_2$ such that for any $\eps\in(0,1),\, l,\d>0$
and any cylinder
\[
Q_\rho^{(\d)}(y,s)=B_\rho(y)\times
(s-\d^{2-p}\rho^p,\, s+\d^{2-p}\rho^p)\subset \O_T, \quad\rho\le R
\]
and any  $\xi\in C_0^\infty(Q_\rho^{(\d)}(y,s))$ such that $\xi(x,t)=1$ for $(x,t)\in Q_{\rho/2}^{(\d)}(y,s)$
\begin{eqnarray}
\nonumber
&&I_1:=\d^2 \int_{L(t)}u(x,t)G\left(\frac{u(x,t)-l}{\d}\right)\xi(x,t)^kdx\\
\nonumber
&+&\iint_{L}\left\{\int^u_l\left(1+\frac{s-l}{\d}\right)^{-1+\l}\left(\frac{s-l}{\d}\right)^{-2\l}ds+
u\left(1+\frac{u-l}{\d}\right)^{-1+\l}\left(\frac{u-l}{\d}\right)^{-2\l}\right\}|\n u|^p\xi(x,\tau)^kdx\,d\tau\\
\nonumber
&\le& \g \d^2 \iint_L u \left(\frac{u-l}{\d}\right)|\xi_t|\xi^{k-1}dxd\tau
+\g \d^p\iint_L u\left[\left(1+\frac{u-l}{\d}\right)^{1-\l}\left(\frac{u-l}{\d}\right)^{2\l}\right]^{p-1}|\n \xi|^p \xi^{k-p}dxd\tau
\\
\nonumber
&+&\eps\d^p \iint_L u\left(1+\frac{u-l}{\d}\right)^{(1+\l)(p-1)}|\n\xi|^p \xi^{k-p}dxd\tau
+\g \iint_L u^p \int_l^u\left(1+\frac{s-l}{\d}\right)^{-1+\l}\left(\frac{s-l}{\d}\right)^{-2\l}ds\,F_2(x) \xi^k dxd\tau
\\
\nonumber
&+& \g \eps^{-\frac1{p-1}}\iint_{L} u^{p+1}\left(1+\frac{u-l}{\d}\right)^{-1-\l}F_1(x)\xi^kdxd\tau
+
\g\eps^{-\frac1{p-1}}\rho^p\d^{2-p}\int_{B_\rho(y)} f_1(x)^\frac{p}{p-1} dx\\
&&\ \ \ \ \ \ \ \ \ \ \ \ \  \  +\g \rho^p\d^{3-p} \int_{B_\rho(y)}
f_2(x) dx:=\sum_{i=2}^8 I_i, \label{e2.3b}
\end{eqnarray}
where $L=Q_\rho^{(\d)}(y,s)\cap\{u>l\}$, $L(t)=L\cap \{\tau=t\}$ and $\l=\min\{\frac1{pn}, \frac{p-2}2\}$, $k=\frac{(p+2)(p-1)(1+\l)}{p-1-\l}+p$.
\end{lemma}
\proof First, note that
\begin{equation}
\label{e2.4b}
\int_l^u\left(1+\frac{s-l}{\d}\right)^{-1+\l}\left(\frac{s-l}{\d}\right)^{-2\l}ds\le \g \d,
\end{equation}
and
\begin{eqnarray}
\nonumber
\int_l^u wdw\int_l^w\left(1+\frac{s-l}{\d}\right)^{-1+\l}\left(\frac{s-l}{\d}\right)^{-2\l}ds
\nonumber
=\frac12\int_l^u\left(1+\frac{s-l}{\d}\right)^{-1+\l}\left(\frac{s-l}{\d}\right)^{-2\l}(u^2-s^2)ds\\
\nonumber
\ge\frac14 u(u-l)\int_l^{\frac{u+l}2}\left(1+\frac{s-l}{\d}\right)^{-1+\l}\left(\frac{s-l}{\d}\right)^{-2\l}ds
=\frac{\d^2}4u\left(\frac{u-l}{\d}\right)\int_0^{\frac{u-l}{2\d}}(1+z)^{-1+\l}z^{-2\l}dz\\
\label{e2.5b}
\ge \g\d^2uG\left(\frac{u-l}{\d}\right).
\end{eqnarray}
Test \eqref{e1.12b} by $\var$ defined by
\begin{equation}
\var(x,t)=u(x,t)\left[\int_l^{u(x,t)}\left(1+\frac{s-l}{\d}\right)^{-1+\l}\left(\frac{s-l}{\d}\right)^{-2\l}ds\right]_+\xi(x,t)^k,
\end{equation}
and $t_1=s-\d^{2-p}\rho^p$, $t_2=t$.
Using \eqref{e1.2b} we have for any $t>0$
\begin{eqnarray*}
&&\int_{L(t)}\int_l^uwdw\int_l^w\left(1+\frac{s-l}{\d}\right)^{-1+\l}\left(\frac{s-l}{\d}\right)^{-2\l}ds \xi^k dx\\
&+&\iint_L\left\{\int_l^u\left(1+\frac{s-l}{\d}\right)^{-1+\l}\left(\frac{s-l}{\d}\right)^{-2\l}ds
+u\left(1+\frac{u-l}{\d}\right)^{-1+\l}\left(\frac{u-l}{\d}\right)^{-2\l}\right\}|\n u|^p\xi^k dxdt\\
&\le& \g \iint_L \int_l^uwdw \int_l^w\left(1+\frac{s-l}{\d}\right)^{-1+\l}\left(\frac{s-l}{\d}\right)^{-2\l}ds|\xi_t|\xi^{k-1}dxdt\\
&+&\g \iint_Lu\int_l^u\left(1+\frac{s-l}{\d}\right)^{-1+\l}\left(\frac{s-l}{\d}\right)^{-2\l}ds\left[|\n u|^{p-1}+g_1 u^{p-1} +f_1\right]|\n\xi|\xi^{k-1}dxdt\\
&+&\g \iint_Lu\int_l^u\left(1+\frac{s-l}{\d}\right)^{-1+\l}\left(\frac{s-l}{\d}\right)^{-2\l}ds\left[ h |\n u|^{p-1}+g_2 u^{p-1}+f_2\right]\xi^kdxdt.
\end{eqnarray*}
From this using  \eqref{e2.4b}, \eqref{e2.5b} and Young's inequality we obtain the required~\eqref{e2.3b}.~\qed

\bigskip
Set
\begin{equation}
\label{e2.6b}
\psi(x,t)=\frac1{\d}\left[\int_l^{u(x,t)}s^\frac1{p}\left(1+\frac{s-l}{\d}\right)^{-\frac{1-\l}{p}}\left(\frac{s-l}{\d}\right)^{-\frac{2\l}{p}}ds\right]_+.
\end{equation}

\begin{lemma}
\label{lem2.3b}
Let the conditions of Lemma~\ref{lem2.2b} be fulfilled.
Then there exists $\nu_1\in(0,1)$ depending only on $n,p,c_1,c_2$ such that
the inequality
\begin{equation}
\label{e2.8-9b}
\calF_1(2R)+\calF_2(2R)
\le\nu_1
\end{equation}
implies that
\begin{eqnarray}
\nonumber \int_{L(t)}uG\left(\frac{u-l}{\d}\right)\xi^k
dx+\d^{p-2}\iint_L |\n\psi|^p \xi^k dxd\tau
\le\frac{\nu_1^\frac12\d^{p-2}}{\rho^p}\iint_L u\left(1+\frac{u-l}{\d}\right)^{(1+\l)(p-1)}\xi^{k-p}dxd\tau\ \ &&\\
\nonumber
+\g \frac{\d^{p-2}}{\rho^p}\iint_L u\left(1+\frac{u-l}{\d}\right)^{(1-\l)(p-1)}
\left(\frac{u-l}{\d}\right)^{2\l(p-1)}\xi^{k-p}dxd\tau\ \ &&\\
\nonumber
+\g \nu_1^{-\frac12}\frac{l^{p+1}\rho^p}{\d^p}\int_{B_\rho(y)} F_1(x) dx
+\g\frac{l^{p}\rho^p}{\d^{p-1}}\int_{B_\rho(y)} F_2(x) dx\ \ &&\\
\label{e2.10b}
+\g \frac{\nu_1^{-\frac12}\rho^p}{\d^p}\int_{B_\rho(y)}  f_1(x)^\frac{p}{p-1} dx
+\g \frac{\rho^p}{\d^{p-1}}\int_{B_\rho(y)}  f_2(x) dx.\ \ &&
\end{eqnarray}
\end{lemma}
\proof
In the notation of Lemma~\ref{lem2.2b} with $\eps=\nu_1^\frac{p-1}2$,  using the Young inequality we have
\begin{eqnarray}
\nonumber
I_2+I_3+I_4\le \nu_1^\frac12 \frac{\d^p}{\rho^p}\iint_L u \left(1+\frac{u-l}{\d}\right)^{(1+\l)(p-1)}\xi^{k-p}dxd\tau\ \ && \\
\label{e2.12b}
+\g \frac{\d^p}{\rho^p}\iint_L u\left[\left(1+\frac{u-l}{\d}\right)^{1-\l}\left(\frac{u-l}{\d}\right)^{2\l}\right]^{p-1}
\xi^{k-p}dxd\tau.\ \ &&
\end{eqnarray}
Set
\begin{eqnarray*}
I_9&=&\g \iint_L (u-l)_+^p\int_l^u\left(1+\frac{s-l}{\d}\right)^{-1+\l}\left(\frac{s-l}{\d}\right)^{-2\l}ds\,F_2(x)\xi^k dxd\tau,\\
I_{10}&=& \g \nu_1^{-\frac12}\iint_L (u-l)_+^{p+1}\left(1+\frac{u-l}{\d}\right)^{-1-\l}F_1(x)\xi^k dxd\tau.
\end{eqnarray*}
First we estimate $I_{10}$. By the Young inequality and Lemma~\ref{lem2.1b} we obtain
\begin{eqnarray}
\nonumber
I_{10} \le \g \nu_1^\frac12\iint_L (u-l)_+\left(1+\frac{u-l}{\d}\right)^{-1-\l}|\n u |^p \xi^k dxd\tau\\
\nonumber
+\g \nu_1^\frac12 \d^{-p}\iint_L(u-l)_+^{p+1}\left(1+\frac{u-l}{\d}\right)^{-1-\l-p}|\n u |^p \xi^k dxd\tau\\
\nonumber
+\g \frac{\nu_1^\frac12}{\rho^p}\iint_L (u-l)_+^{p+1}\left(1+\frac{u-l}{\d}\right)^{-1-\l}\xi^{k-p}dxd\tau\\
\le \frac14 I_1 +\g \nu_1^\frac12 \frac{\d^p}{\rho^p}\iint_L u \left(1+\frac{u-l}{\d}\right)^{(1+\l)(p-1)}\xi^{k-p}dxd\tau.\label{e2.13b}
\end{eqnarray}

To estimate $I_9$ we consider the weak solution to the problem
\[
-\Delta_p H =F_2, \quad H\in \W^{1,p}(B_\rho(y)),
\]
i.e.
\begin{equation}
\label{e2.14b}
\int_{B_\rho(y)}|\n H|^{p-2}\n H\cdot \n \var dx = \int_{B_\rho(y)}F_2\var dx,
\end{equation}
for any $\var\in \W^{1,p}(B_\rho(y))$. By \cite{KiMa} we have
\begin{equation}
\label{e2.15b}
\| H \|_{L^\infty(B_\rho(y))}\le \g\, \calF_2 (2\rho).
\end{equation}
Testing \eqref{e2.14b} by
\[
\var=(u-l)_+^p\int_l^u\left(1+\frac{s-l}{\d}\right)^{-1+\l}\left(\frac{s-l}{\d}\right)^{-2\l}ds\xi^k,
\]
and using the Young inequality we have
\begin{eqnarray}
\nonumber
\int_{B_\rho(y)}(u-l)_+^p\int_l^u\left(1+\frac{s-l}{\d}\right)^{-1+\l}\left(\frac{s-l}{\d}\right)^{-2\l}ds\, F_2(x)\xi^k dx\\
\nonumber
\le \g\int_{B_\rho(y)}|\n H|^{p-1}\left|\n\left\{(u-l)^p\int_l^u\left(1+\frac{s-l}{\d}\right)^{-1+\l}\left(\frac{s-l}{\d}\right)^{-2\l}ds
\xi^k\right\}\right|dx
\end{eqnarray}
\begin{eqnarray}
\nonumber
\le \frac{\nu_1^\frac12}{8\g}\int_{B_\rho(y)}\int_l^u\left(1+\frac{s-l}{\d}\right)^{-1+\l}\left(\frac{s-l}{\d}\right)^{-2\l}ds
|\n u|^p \xi^k dx\\
\nonumber
+\frac{\nu_1^\frac12}{8\g}\int_{B_\rho(y)}\int_l^u(u-l)_+\left(1+\frac{s-l}{\d}\right)^{-1+\l}\left(\frac{s-l}{\d}\right)^{-2\l}ds
|\n u|^p \xi^k dx\\
\nonumber
+\g \nu_1^{-\frac1{2(p-1)}}\int_{B_\rho(y)}(u-l)_+^p\int_l^u\left(1+\frac{s-l}{\d}\right)^{-1+\l}\left(\frac{s-l}{\d}\right)^{-2\l}ds
|\n H|^p \xi^k dx\\
\nonumber
+\g \nu_1^{-\frac1{2(p-1)}}\int_{B_\rho(y)}(u-l)_+^{p+1}
\left(1+\frac{u-l}{\d}\right)^{-1+\l}\left(\frac{u-l}{\d}\right)^{-2\l}|\n H|^p \xi^k dx\\
\label{e2.16b}
+\frac{\nu_1\d}{\rho^p}\int_{B_\rho(y)}(u-l)_+^p\xi^{k-p}dx.
\end{eqnarray}
Using the definition of the weak solution to \eqref{e2.14b} again, we have
\begin{eqnarray}
\nonumber
\int_{B_\rho(y)}(u-l)_+^p\int_l^u\left(1+\frac{s-l}{\d}\right)^{-1+\l}\left(\frac{s-l}{\d}\right)^{-2\l}ds |\n H |^p \xi^k dx\\
\nonumber
=\int_{B_\rho(y)}F_2(x)H(u-l)_+^p\int_l^u\left(1+\frac{s-l}{\d}\right)^{-1+\l}\left(\frac{s-l}{\d}\right)^{-2\l}ds  \xi^k dx\\
\label{e2.17b}
+\int_{B_\rho(y)}H |\n H |^{p-2}\n H \cdot \n\left\{(u-l)_+^p\int_l^u\left(1+\frac{s-l}{\d}\right)^{-1+\l}\left(\frac{s-l}{\d}\right)^{-2\l}ds
\xi^k \right\}dx,
\end{eqnarray}
and
\begin{eqnarray}
\nonumber
\int_{B_\rho(y)}(u-l)_+^{p+1}
\left(1+\frac{u-l}{\d}\right)^{-1+\l}\left(\frac{u-l}{\d}\right)^{-2\l}|\n H|^p \xi^k dx\\
\nonumber
=\int_{B_\rho(y)}H(u-l)_+^{p+1}
\left(1+\frac{u-l}{\d}\right)^{-1+\l}\left(\frac{u-l}{\d}\right)^{-2\l} F_2(x)\xi^k dx\\
\label{e2.18b}
+\int_{B_\rho(y)}H |\n H |^{p-2}\n H \cdot \n\left\{(u-l)_+^{p+1}
\left(1+\frac{u-l}{\d}\right)^{-1+\l}\left(\frac{u-l}{\d}\right)^{-2\l}\xi^k\right\}dx.
\end{eqnarray}

The terms in the right hand side of \eqref{e2.17b} have been estimated in \eqref{e2.16b}.
The right hand side of \eqref{e2.18b} is estimated similarly to \eqref{e2.13b} using the Young inequality.
Thus using \eqref{e2.8-9b} and \eqref{e2.15b} and collecting \eqref{e2.12b},\eqref{e2.13b}, \eqref{e2.16b}-\eqref{e2.18b}
we arrive at the required \eqref{e2.10b}.\qed

\section{Local boundedness of solutions. Proof of Theorem~\ref{thm1.1b}}
\label{bdd}

Let $(x_0,t_0)$ be an arbitrary point in $\O_T$.
Let
\[
R\le \frac12\min\left\{1,{\rm dist}\,(x_0,\partial\O_T),
t_0^\frac12, (T-t_0)^\frac12\right\}.
\]
Let
\[
Q_R(x_0,t_0)=B_R(x_0)\times (t_0-R^2,t_0+R^2).
\]
Fix a point $(y,s)\in Q_\frac{R}2(x_0,t_0)$. For $j=1,2,\dots$ set
\[
\rho_j=R 2^{-j},\quad Q_j=B_j\times (s-\d_j^{2-p}\rho_j^p,s+\d_j^{2-p}\rho_j^p),
\quad B_j=B_{\rho_j}(y),\quad L_j=Q_j\cap\O_T\cap\{u(x,t)>l_j\}.
\]
Let $\xi_j\in C_0^\infty(Q_j)$ be such that $\xi_j(x,t)=1$ for $(x,t)\in B_{j+1}\times (s-\frac34 \d_j^{2-p}\rho_j^p,\,s+\frac34 \d_j^{2-p}\rho_j^p)$,
$|\n \xi_j|\le \g \rho_j^{-1}$, $|\frac{\partial \xi_j}{\partial t}|\le \g \d_j^{p-2}\rho_j^{-p}$.

The sequences of positive numbers $(l_j)_{j\in\N}$ and $(\d_j)_{j\in\N}$ are defined inductively as follows.

Set $l_0=1$ and assume that $l_1,l_2,\dots,l_j$ and $\d_0,\d_1,\dots,\d_{j-1}$ have been already chosen.
Let us show  how to chose $l_{j+1}$ and $\d_j$.

Define the sequence $(\a_j)_{j\in\N}$ by
\begin{equation}
\label{e3.2b}
\a_j=
\rho_j
+
\left(\frac1{\rho_j^{n-p}}\int_{B_j}F_1(x)dx\right)^\frac1{p}+
\left(\frac1{\rho_j^{n-p}}\int_{B_j}F_2(x)dx\right)^\frac1{p-1}.
\end{equation}

For $l\ge l_j+\a_j$ set
\begin{equation}
\label{e3.1b}
A_j(l)=\frac{(l-l_j)^{p-2}}{\rho_j^{n+p}}\iint_{\widetilde L_j} \frac{u}{l}\left(\frac{u-l_j}{l-l_j}\right)^{(1+\l)(p-1)}\xi_j^{k-p}dxd\tau
+\sup_t \frac1{\rho_j^n}\int_{\widetilde L_j(t)}\frac{u}{l}G\left(\frac{u-l_j}{l-l_j}\right)\xi_j^k dx,
\end{equation}
where $\widetilde L_j=\widetilde Q_j\cap \O_T\cap\{u(x,t)>l_j\}$,
$\widetilde Q_j=B_j\times (s-(l-l_j)^{2-p}\rho_j^p,s+(l-l_j)^{2-p}\rho_j^p)$.


Fix a positive number $\varkappa\in (0,1)$ depending on $n,p,c_1,c_2$, which will be specified later.
If
\begin{equation}
\label{e3.3b}
A_j(l_j+\a_j)\le \varkappa,
\end{equation}
we set $l_{j+1}=l_j+\a_j$.

Note that $A_j(l)\searrow 0$ as $l\to \infty$.
So
if
\begin{equation}
\label{e3.4b}
A_j(l_j+\a_j)> \varkappa,
\end{equation}
there exists $\bar l>l_j+\a_j$ such that $A_j(\bar l)=\varkappa$. In this case we set $l_{j+1}=\bar l$.

In both cases we set $\d_j=l_{j+1}-l_j$. Note that our choices guarantee that $\widetilde Q_j\subset Q_R(x_0,t_0)$ and
\begin{equation}
\label{e3.5b}
A_j(l_{j+1})\le \varkappa.
\end{equation}

\begin{lemma}
\label{lem3.1b}
Let the conditions of Theorem~\ref{thm1.1b} be fulfilled. The for all $j\ge 1$
there exists $\g>0$ depending on the data, such that
\begin{equation}
\label{e3.6b}
\d_j\le \frac12 \d_{j-1}+\a_j+\g (1+l_j)\nu_1^{-\frac1{2p}}\left(\rho_j^{p-n}\int_{B_j}F_1(x)dx\right)^\frac1{p}
+\g (1+l_j)\left(\rho_j^{p-n}\int_{B_j}F_2(x)dx\right)^\frac1{p-1}.
\end{equation}
\end{lemma}
\proof Fix $j\ge 1$. Without loss assume that
\begin{equation}
\label{e3.7b}
\d_j>\frac12 \d_{j-1},\quad \d_j>\a_j,
\end{equation}
since otherwise \eqref{e3.6b} is evident. The second inequality in \eqref{e3.7b} guarantees that $A_j(l_{j+1})=\varkappa$ and $\tilde Q_j=Q_j$.

Let us estimate the terms in the right hand side of \eqref{e3.1b} with $l=l_{j+1}$. For this we decompose $L_j$ as
$L_j=L_j\pr\cup L_j^{\prime\prime}$,
\[
L_j\pr=\left\{x\in L_j\,:\,\frac{u(x)-l_j}{\d_j}<\eps_1\right\}, \quad L_j^{\prime\prime}=L_j\setminus L\pr_j,
\]
where $\eps_1$ depending on $n,p,c_1,c_2$ is small enough to be determined later.

We also have
\begin{equation}
\label{e3.8b}
l_j\ge 1.
\end{equation}
Recall that $\xi_{j-1}=1$ on $Q_j$. By \eqref{e3.5b} we have
\begin{eqnarray}
\nonumber
\frac{\d_j^{p-2}}{\rho_j^{n+p}}\iint_{L\pr_j}\frac{u}{l_{j+1}}\left(\frac{u-l_j}{\d_j}\right)^{(1+\l)(p-1)}\xi_j^{k-p}dxd\tau
\le 2\eps_1^{(1+\l)(p-1)}\rho_j^{-n}\sup_t \int_{L\pr_j(t)}\frac{u}{l_{j+1}}\xi_{j-1}^k dx\\
\label{e3.9b}
\le 2^{n+1} \eps_1^{(1+\l)(p-1)}\rho_{j-1}^{-n}\sup_t\int_{L_{j-1}(t)}\frac{u}{l_j}G\left(\frac{u-l_{j-1}}{\d_{j-1}}\right)\xi_{j-1}^k dx
\le 2^{n+1} \eps_1^{(1+\l)(p-1)} \varkappa.
\end{eqnarray}

Let
\[
\psi_j(x,t)=\frac1{\d_j}\left(\int_{l_j}^{u(x,t)}s^\frac1{p}\left(1+\frac{s-l_j}{\d_j}\right)^{-\frac{1-\l}{p}}
\left(\frac{s-l_j}{\d_j}\right)^{-\frac{2\l}{p}}ds\right)_+.
\]
Then by the Young inequality
\begin{eqnarray}
\nonumber
&&\iint_{L_j^{\prime\prime}}u\left(\frac{u-l_j}{\d_j}\right)^{(1+\l)(p-1)}\xi_j^{k-p}dxd\tau
\le \eps_2\iint_{L_j^{\prime\prime}}u(x,\tau)dxd\tau\\
\label{e3.10b}
&+& \g(\eps_2) \iint_{L_j^{\prime\prime}}u\left(\frac{u-l_j}{\d_j}\right)^{(1+\l)(p-1)z}\xi_j^{(k-p)z}dxd\tau,
\end{eqnarray}
where
\[
z=\frac{n+\rho(\l)}{n}\frac{p-1-\l}{(1+\l)(p-1)}, \quad \rho(\l)=\frac{p}{p-1-\l}.
\]
Similarly to \eqref{e3.9b} we have
\begin{equation}
\label{e3.11b}
\iint_{L_j^{\prime\prime}}u(x,\tau)dxd\tau \le 2^{n+1} \frac{\rho_j^{p+n}}{\d_j^{p-2}}l_j\varkappa.
\end{equation}
Using the evident inequality
\[
c(\eps_1)^{-1}\psi_j(x,t)^{\rho(\l)}\le u(x,t)^\frac1{p-1-\l}\left(\frac{u(x,t)-l_j}{\d_j}\right) \le c(\eps_1)\psi_j(x,t)^{\rho(\l)},
\quad (x,t)\in L_j^{\prime\prime},
\]
the Sobolev inequality  and Lemma~\ref{lem2.3b} with $l=l_j$, $\d=\d_j$, we obtain
\begin{eqnarray}
\nonumber
 \iint_{L_j^{\prime\prime}}u\left(\frac{u-l_j}{\d_j}\right)^{(1+\l)(p-1)z}\xi_j^{(k-p)z}dxd\tau
 \le \g(\eps_1) \iint_{L_j^{\prime\prime}}u^{-\frac{\rho(\l)}{n}}\psi_j^{p\frac{n+\rho(\l)}{n}}\xi_j^{(k-p)z}dxd\tau\\
 \nonumber
 \le \g(\eps_1) l_j^{-\frac{\rho(\l)}{n}}\iint_{L_j^{\prime\prime}}\psi_j^{p\frac{n+\rho(\l)}{n}}\xi_j^{(k-p)z}dxd\tau\\
 \nonumber
  \le \g(\eps_1) l_j^{-\frac{p}{n}}\left(\sup_t \int_{L_j^{\prime\prime}(t)}u^\frac{p-2-\l}{p-1-\l}\psi_j^{\rho(\l)}\xi_j^\frac{(k-p)zn}{n+\rho(\l)}dx\right)^\frac{p}{n}
 \iint_{L_j}\left|\n\left(\psi_j \xi_j^\frac{(k-p)zn}{(n+\rho(\l))p}\right)\right|^pdxd\tau
 \end{eqnarray}
 \begin{eqnarray}
 \nonumber
 \le \g(\eps_1) \d_j^{2-p}l_j^{-\frac{p}{n}}
 \left[
 \frac{\d_j^{p-2}}{\rho_j^p}\iint_{L_j} u \left(1+\frac{u-l_j}{\d_j}\right)^{(1+\l)(p-1)}\xi_jdxd\tau
\right.
 \\
 \nonumber
 +
 \g \frac{\d_j^{p-2}}{\rho_j^p}
 \iint_{L_j}
 u\left(1+\frac{u-l_j}{\d_j}\right)^{(1-\l)(p-1)}
\left(\frac{u-l_j}{\d_j}\right)^{2\l(p-1)} \xi_jdxd\tau\\
\nonumber
+\nu_1^{-\frac12}l_j^{p+1}
\left(\frac{\rho_j}{\d_j}\right)^p
\int_{B_j}F_1dx
+l_j^{p}{\d_j}\left(\frac{\rho_j}{\d_j}\right)^p\int_{B_j}F_2dx
\\
\left.
\label{e3.12b}
+{\nu_1^{-\frac12}}\left(\frac{\rho_j}{\d_j}\right)^p\int_{B_j}f_1^\frac{p}{p-1}dx
+{\d_j}\left(\frac{\rho_j}{\d_j}\right)^p\int_{B_j}f_2dx\right]^{1+\frac{p}{n}}.
\end{eqnarray}
From \eqref{e3.9b}-\eqref{e3.12b} and from the fact that $\xi_{j-1}=1$ on $Q_{j-1}$, we obtain
\begin{eqnarray}
\nonumber
\frac{\d_j^{p-2}}{\rho_j^{n+p}}\iint_{L_j}\frac{u}{l_{j+1}}\left(\frac{u-l_j}{\d_j}\right)^{(1+\l)(p-1)}\xi_j^{k-p}dxd\tau
\le (2^{n+1}
\eps_1^{(1+\l)(p-1)}+2^{n+1}\eps_2)\varkappa
\\
\nonumber
+
\g(\eps_1,\eps_2)\left[\varkappa+\nu_1^{-\frac12}\frac{l_j^{p}\rho_j^{p-n}}{\d_j^p}\int_{B_j} F_1(x)dx
+\frac{l_j^{p-1}\rho_j^{p-n}}{\d_j^{p-1}}\int_{B_j} F_2(x)dx
\right.
\\
\label{e3.13b}
\left.
+\nu_1^{-\frac12}\frac{\rho_j^{p-n}}{\d_j^p}\int_{B_j} f_1(x)^\frac{p}{p-1}dx +
\frac{\rho_j^{p-n}}{\d_j^{p-1}}\int_{B_j}f_2(x)dx\right]^{1+\frac{p}{n}}.
\end{eqnarray}

Let us estimate the second term in the right hand side of \eqref{e3.1b}.
By Lemma~\ref{lem2.3b} with $l=l_j, \ \d=\d_j$ we have
\begin{eqnarray}
\nonumber
\sup_t\int_{L_j(t)} u G\left(\frac{u-l_j}{\d_j}\right)\xi_j^k dx
\le \g \left(\nu_1+
\eps_1^{2\l(p-1)}\right)
\d_j^{p-2}\rho_j^{-p}\iint_{L\pr_j}u\xi_j^{k-p}dxd\tau\\
\nonumber
+\g(\eps_1)\d_j^{p-2}\rho_j^{-p}
\iint_{L_j^{\prime\prime}}u\left(\frac{u-l_j}{\d_j}\right)^{(1+\l)(p-1)}\xi_j^{k-p}dxd\tau
+
\nu_1^{-\frac12}\frac{l_j^{p}\rho_j^p}{\d_j^p}\int_{B_j} F_1(x)dx + \g \frac{l_j^{p-1}\rho_j^p}{\d_j^{p-1}}\int_{B_j} F_2(x)dx\\
\label{e3.14b}
+\nu_1^{-\frac12}\frac{\rho_j^p}{\d_j^p}\int_{B_j}f_1(x)^\frac{p}{p-1}dx +\g\frac{\rho_j^p}{\d_j^{p-1}}\int_{B_j} f_2(x)dx.
\end{eqnarray}
The first two terms of the right hand side of \eqref{e3.14b} were estimated in \eqref{e3.9b} and in \eqref{e3.12b}.
Therefore we conclude from \eqref{e3.13b}, \eqref{e3.14b} that

\begin{eqnarray}
\nonumber
&&\varkappa\le \left(2^{n+1} \eps_1^{(1+\l)(p-1)}+2^{n+1}\eps_2+\nu_1+\eps_1^{2\l(p-1)}\right)\varkappa
+ \g(\eps_1,\eps_2)\left[
\nu_1^{-\frac12}\frac{l_j^{p}\rho_j^{p-n}}{\d_j^p}\int_{B_j} F_1(x)dx\right.\\
\nonumber
&&\left.
+\frac{l_j^{p-1}\rho_j^{p-n}}{\d_j^{p-1}}\int_{B_j}F_2(x)dx
+\nu_1^{-\frac12}\frac{\rho_j^{p-n}}{\d_j^p}\int_{B_j} f_1(x)^\frac{p}{p-1}dx +\frac{\rho_j^{p-n}}{\d_j^{p-1}}\int_{B_j} f_2(x)dx\right]\\
\nonumber
&&+ \g(\eps_1,\eps_2)\left[ \varkappa +\nu_1^{-\frac12}\frac{l_j^{p}\rho_j^{p-n}}{\d_j^p}\int_{B_j} F_1(x)dx
\label{e3.15b}
+\frac{l_j^{p-1}\rho_j^{p-n}}{\d_j^{p-1}}\int_{B_j}F_2(x)dx\right.\\
&&\left.+\nu_1^{-\frac12}\frac{\rho_j^{p-n}}{\d_j^p}\int_{B_j} f_1(x)^\frac{p}{p-1}dx +\frac{\rho_j^{p-n}}{\d_j^{p-1}}\int_{B_j} f_2(x)dx\right]^{1+\frac{p}{n}}.
\end{eqnarray}

Choose $\nu_1<\frac\varkappa{16}$, $\eps_1,\eps_2,$ such that
\[
2^{n+1} \eps_1^{(1+\l)(p-1)}+2^{n+1}\eps_2+\eps_1^{2\l(p-1)}<\frac1{16},
\]
and $\varkappa$ such that $\g(\eps_1,\eps_2)\varkappa^\frac{p}{n}=\frac1{16}$.
%
Hence \eqref{e3.15b} yields \eqref{e3.6b} which completes the proof of the lemma.\qed

In order to complete the proof of Theorem~\ref{thm1.1b} we sum up \eqref{e3.6b} with respect to $j$ from 1 to $J-1$
\begin{eqnarray}
\nonumber
&&l_J\le \g \d_0 +\g \sum_{j=1}^\infty \a_j +\g(1+l_J)\nu_1^{-\frac1{2p}}\sum_{j=1}^\infty \left(\rho_j^{p-n}\int_{B_j}F_1(x)dx\right)^\frac1{p}
+\g(1+l_J)\sum_{j=1}^\infty\left(\rho_j^{p-n}\int_{B_j}F_2(x)dx\right)^\frac1{p-1}\\
\nonumber
&&\le  \g \d_0+\g R + \g  (1+l_J)\left[
\nu_1^{-\frac1{2p}}\int_0^{2R}\left(\frac1{r^{n-p}}\int_{B_r(y)}F_1(x)dx\right)^\frac1{p}\frac{dr}r\right.\\
\label{e3.16b}
 &&\left.\ \ \ \ \ \ \ +\int_0^{2R}\left(\frac1{r^{n-p}}\int_{B_r(y)}F_2(x)dx\right)^\frac1{p-1}\frac{dr}r\right].
\end{eqnarray}

Let us estimate $\d_0$.
If $l_1=l_0+\a_0$ then
$\d_0=\a_0$.
If on the other hand $l_1,\d_0$ are defined by $A(l_1)=\varkappa$ then by \eqref{e3.1b}
\begin{equation}
\label{e3.17b}
\d_0\le \g  \left(\frac1{R^{p+n}}\iint_{Q_R(x_0,t_0)} u^{p+\l(p-1)}dxdt\right)^\frac{1}{2+\l(p-1)}
+ \g \sup_{t_0-R^2<t<t_0+R^2} \left(\frac1{R^{n}}\int_{B_R(x_0)}u^2dx\right)^\frac{1}{2}.
\end{equation}

Now we choose $R_0>0$ such that $\calF_1(R_0)+\calF_2(R_0)<\nu_1\le \frac12\g^{-1}\nu_1^\frac1{2p}$, where $\g$ as in last line of \eqref{e3.16b}.
Then by \eqref{e3.16b}, for $R\le R_0$ we obtain
\begin{eqnarray}
\nonumber
l_J\le \g  +\g \left(\frac1{R^{p+n}}\iint_{Q_R(x_0,t_0)} u^{p+\l(p-1)}dxdt\right)^\frac{1}{2+\l(p-1)}
+ \g \sup_{t_0-R^2<t<t_0+R^2} \left(\frac1{R^{n}}\int_{B_R(x_0)}u^2dx\right)^\frac{1}{2}\\
\label{e3.18b}
+\g R
+\g \int_0^{2R}\left(\frac1{r^{n-p}}\int_{B_r(y)}F_1(x)dx\right)^\frac1{p}\frac{dr}r
+\g \int_0^{2R}\left(\frac1{r^{n-p}}\int_{B_r(y)}F_2(x)dx\right)^\frac1{p-1}\frac{dr}r.
\end{eqnarray}
Hence the sequence $(l_j)_{j\in\N}$ is convergent, and $\d_j\to 0\,(j\to\infty)$, and we can pass to the limit $J\to\infty$ in \eqref{e3.18b}.
Let $l=\lim_{j\to\infty}l_j$. From \eqref{e3.5b} we conclude that
\begin{equation}
\label{e3.19b}
\frac1{\rho_j^{n+p}}\iint_{Q_j}u(u-l)_+^{(1+\l)(p-1)}\le \g \vark\, l \d_j^{1+\l(p-1)} \to 0\quad (j\to \infty).
\end{equation}
Choosing $(y,s)$ as a Lebesgue point of the function $u(u-l)_+^{(1+\l)(p-1)}$ we conclude that $u(y,s)\le l$ and hence $u(y,s)$
is estimated from above by the right hand side of \eqref{e3.18b}. Applicability of the Lebesgue differentiation theorem
follows from \cite[Chap.~II,\,Sec.~3]{Guz}.

Taking essential supremum over $Q_{R/2}(x_0,t_0)$ we complete the proof.\qed
\section{Proof of Theorem~\ref{thm1.4}}
\label{degiorgi}
In this section we prove the
Theorem~\ref{thm1.4} which is
a DeGiorgi-type lemma \cite{DiGV}.
Here we assume the structure conditions
\begin{eqnarray}
\nonumber
\A (x,t,u, \zeta)\zeta &\ge& c_1 |\zeta|^p,\quad \zeta\in\R^n,\\
|\A (x,t,u, \zeta)|&\le& c_2|\zeta|^{p-1}+f_1(x),\label{e1.2}\\
\nonumber |a_0(x,t,u,\zeta)|&\le& |\zeta|^{p-1}+f_2(x),
\end{eqnarray}
with some positive constants $c_1, c_2$ and nonnegative functions
$f_1(x),f_2(x)$. These assumptions follow from \eqref{e1.2b} due to the boundedness of $u$ and $h=1$.

We assume that
\[
f_1^\frac{p}{p-1}\in \widetilde K_p,\quad f_2\in K_p.
\]

We provide the proof of \eqref{e1.13}, while the proof of \eqref{e1.15} is completely similar.

Set $v=u-\mu_-$, $M=\esssup_{\O_T} |u(x,t)|$. In the sequel $\g$
will denote a constant depending on the data and $M$, which, as usual, can vary from line to line.

\medskip

\begin{lemma}
\label{lem3.1}
Let $u$ be a solution to \eqref{e1.1}. Then for any $l,\d>0$ and $\eps\in(0,1)$
and any cylinder
\[
Q_\rho^{(\d)}(y,s)=B_\rho(y)\times
(s-\d^{2-p}\rho^p,s+\d^{2-p}\rho^p)\subset \O_T, \quad\rho\le R
\]
and any  $\xi\in C_0^\infty(Q_\rho^{(\d)}(y,s))$ such that $\xi(x,t)=1$ for $(x,t)\in Q_{\rho/2}^{(\d)}(y,s)$
with $|\n\xi|\le\g\frac1{\rho},\,|\xi_t|\le \g \frac{\d^{p-2}}{\rho^p}$
we have
\begin{eqnarray}
\nonumber
&&\d^2\int_{L(t)}G\left(\frac{l-v(x,t)}{\d}\right)\xi(x,t)^kdx\\
\nonumber
&&+\iint_{L}
\left(1+\frac{l-v}{\d}\right)^{-1+\l}\left(\frac{l-v}{\d}\right)^{-2\l}
|\n v|^p\xi(x,\tau)^kdx\,d\tau\\
\nonumber
&\le& \g \frac{\d^{p}}{\rho^p}\iint_L \left(1+\frac{l-v}{\d}\right)^{(1-\l)(p-1)}\left(\frac{l-v}{\d}\right)^{2\l(p-1)} \xi^{k-p}dxd\tau\\
\nonumber
&+&\eps \frac{\d^{p}}{\rho^p}\iint_L \left(1+\frac{l-v}{\d}\right)^{(1+\l)(p-1)} \xi^{k-p}dxd\tau\\
\label{e3.1}
&+&\g\eps^{-\frac1{p-1}}{\rho^p}{\d^{2-p}}\int_{B_\rho(y)} F_1(x) dx
 +
\g {\rho^p}{\d^{3-p}} \int_{B_\rho(y)} F_2(x)dx,
\end{eqnarray}
where $L=Q_\rho^{(\d)}(y,s)\cap\{v<l\}$, $L(t)=L\cap \{\tau=t\}$ and $\l=\min\{\frac1{p n},\frac{p-2}2\}$, $k=p+\frac{(p+2)(p-1)(1+\l)}{p-1-\l}$ and
$G(v)$ is defined in the previous section.
\end{lemma}
\proof The proof is similar to that of Lemma~\ref{lem2.2b} with the choice of the test function
\begin{equation*}
\var(x,t)=\left[\int^l_{v(x,t)}\left(1+\frac{l-s}{\d}\right)^{-1+\l}\left(\frac{l-s}{\d}\right)^{-2\l}ds\right]_+\xi(x,t)^k.\qed
\end{equation*}

\bigskip
Set
\begin{equation}
\label{e3.6}
w(x,t)=\left(\frac1{\d}\int_{v(x,t)}^l\left(1+\frac{l-s}{\d}\right)^{-\frac1{p}+\frac{\l}{p}}\left(\frac{l-s}{\d}\right)^{-\frac{2\l}{p}}ds\right)_+.
\end{equation}

Note that
\begin{equation}
\label{e3.7}
c(\eps)w(x,t)^{\rho(\l)}\le \frac{l-v}{\d}\le C(\eps)w(x,t)^{\rho(\l)}\quad \text{if}\ \ \frac{l-v}{\d}\ge \eps,
\end{equation}
with
\begin{equation}
\label{e3.8}
\rho(\l)=\frac{p}{p-1-\l}.
\end{equation}

The next lemma follows from Lemma~\ref{lem3.1} with $\eps=\nu_1$ via the arguments similar to \eqref{e2.13b}--\eqref{e2.18b}.

\begin{lemma}
\label{lem3.2}
Let $u$ be a solution to \eqref{e1.1}. Then for any $l,\d>0$ and $\eps\in(0,1)$
and any cylinder
\[
Q_\rho^{(\d)}(y,s)=B_\rho(y)\times
(s-\d^{2-p}\rho^p,s+\d^{2-p}\rho^p)\subset \O_T, \quad\rho\le R
\]
and any  $\xi\in C_0^\infty(Q_\rho^{(\d)}(y,s))$ such that $\xi(x,t)=1$ for $(x,t)\in Q_{\rho/2}^{(\d)}(y,s)$
with $|\n\xi|\le\g\frac1{\rho},\,|\xi_t| \le \g \frac{\d^{p-2}}{\rho^p}$
there exists $\nu_1\in(0,1)$ depending only on $n,p,c_1,c_2$ such that
the inequality
\begin{equation}
\label{e2.8-9b-f}
\calF_1(2R)+\calF_2(2R)\le\nu_1
\end{equation}
implies that
\begin{eqnarray}
\nonumber
&&\sup_t\int_{L(t)}G\left(\frac{l-v(x,t)}{\d}\right)\xi(x,t)^kdx
+\d^{p-2}\iint_{L}
|\n w|^p\xi(x,\tau)^kdx\,d\tau\\
\nonumber
&\le& \g \nu_1\frac{\d^{p-2}}{\rho^p}\iint_L \left(1+\frac{l-v}{\d}\right)^{(1+{\l})(p-1)} \xi^{k-p}dxd\tau\\
\nonumber
&+&\g \frac{\d^{p-2}}{\rho^p}\iint_L \left(1+\frac{l-v}{\d}\right)^{(1-{\l})(p-1)} \left(\frac{l-v}{\d}\right)^{2\l(p-1)}\xi^{k-p}dxd\tau\\
&+&\g \frac{\rho^p}{\d^p}\int_{B_\rho(y)} F_1(x) dx
\label{e3.9}
+\g \nu_1^{-\frac1{p-1}}\frac{\rho^p}{\d^{p-1}} \int_{B_\rho(y)} F_2(x)dx,
\end{eqnarray}
where $L,\, \l,\, k$ and $G$ are the same as in Lemma~\ref{lem3.1}.
\end{lemma}
Further on we assume that
\begin{equation}
\label{e4.8} \xi\omega\ge B(\rho+\calF_1(2\rho)+\calF_2(2\rho)).
\end{equation}
Let $(x_1,t_1)\in Q_\rho^\theta(y,s)$. Set
\[
r_j=\frac {\rho_0}{4^j}, \rho_0=\frac{\rho}{C}, B_j=B_{r_j}(x_1),
Q_j(l)=B_j \times \left(t_1-\frac{r_j^p}{(l_j-l)^{p-2}},
t_1+\frac{r_j^p}{(l_j-l)^{p-2}}\right), j=1,2,\dots,
\]
$C\ge 16$ will be fixed later depending only on the known data.

Let $\ind_{B_{j+1}}\le \xi_j(x)\le \ind_{B_j}$,
\[
\ind_{\left(t_1-\frac 49 \frac{r_j^p}{(l_j-l)^{p-2}},t_1+\frac49\frac{r_j^p}{(l_j-l)^{p-2}}\right)}\le
\theta_j(t)\le
\ind_{\left(t_1-\frac{r_j^p}{(l_j-l)^{p-2}},
t_1+\frac{r_j^p}{(l_j-l)^{p-2}}, \right)},\quad
\xi_j(x,t)=\xi_j(x)\theta_j(t).\]

We start with the choice of the sequences $l_j, \delta_j, j=0,1,2, \dots$

Set
\begin{equation*}
A_j(l)=\frac{(l_j-l)^{p-2}}{r_j^{n+p}}\iint\limits_{L_j(l)}\left(\frac{l_j-v}{l_j-l}\right)^{(1+\lambda)(p-1)}\xi_j(x,t)^{k-p}dx
dt+\esssup_{t} \frac
1{r_j^n}\int\limits_{L_j(l,t)}G\left(\frac{l_j-v}{l_j-l}\right)\xi_j(x,t)^kdx,
\end{equation*}
where $L_j (l)=Q_j(l)\cap\Omega_T\cap\{v\le l_j\},\quad
L_j(l,t)=L_j\cap\{\tau=t\}$.
Define the sequence
$(\alpha_j)_{j\in\N}$ by
\[
\alpha_j=
r_j+
\int_0^{r_j}\frac{dr}{r}\left(r^{p-n}\int_{B_r(x_1)} F_1(z) dz\right)^\frac1p
+\int_0^{r_j}\frac{dr}{r}\left(r^{p-n}\int_{B_r(x_1)} F_2(z) dz\right)^\frac1{p-1}, \quad j=-1,0,1,2,\dots.
\]
By the definition of the Kato class $\a_j\downarrow 0$ as $j\to \infty$.
Note that
\begin{equation}
\label{alpha-alpha}
\a_{j-1}-\a_j \ge \g \left[\left(r_j^{p-n}\int_{B_j} F_1(z)dz\right)^\frac1p+ \left(r_j^{p-n}\int_{B_j} F_2(z)dz\right)^\frac1{p-1}\right],
\end{equation}
and also
\begin{equation}
\label{alpha-alpha1}
\a_{j-1}-\a_j \le  3r_j+ \g \left[\left(r_{j-1}^{p-n}\int_{B_{j-1}} F_1(z)dz\right)^\frac1p+ \left(r_{j-1}^{p-n}\int_{B_{j-1}} F_2(z)dz\right)^\frac1{p-1}\right].
\end{equation}

Set $l_0=\xi\omega, \quad \overline{l}=\frac{\xi \omega}{2}+\frac{B\alpha_0}{4}$.
Then
$\xi\omega-\overline{l}= \frac
{\xi\omega}2- \frac{B\alpha_0}4\ge \frac {\xi\omega}4\ge
\frac{B\alpha_0}4$ and moreover from \eqref{e1.12} it follows that
\begin{eqnarray}
\nonumber
&&\frac{(\xi\omega-\overline{l})^{p-2}}{r_0^{n+p}}
\iint\limits_{L_0(\overline{l})}\left(\frac{\xi\omega-v}{\xi\omega-\overline{l}}\right)^{(1+\lambda)(p-1)}\xi(x,t)^{k-p}dxdt
\\
\nonumber &\le& \frac{4^{(1+\lambda)(p-1)}(\xi\omega)^{p-2}}{r_0^{n+p}} |\{(x,t)\in Q_0(\overline{l})\,:\, v(x,t)\le \xi\omega\}|\\
\nonumber &\le&
\frac{4^{(1+\lambda)(p-1)}(\xi\omega)^{p-2}}{r_0^{n+p}}|\{(x,t)\in Q_{2\rho}^\theta(y,s)\,:\, u(x,t)\le \mu_-+\xi\omega\}|\\
&\le& \nu 4^{(1+\lambda)(p-1)}(2C)^{n+p}\theta(\xi\omega)^{p-2}.
 \label{e4.9}
\end{eqnarray}
By Lemma \ref{lem3.2}
\begin{eqnarray}
\nonumber
&&  \esssup_{t}\frac{1}{r_0^n}\int\limits_{L_0(\overline{l},t)}G\left(\frac{\xi\omega-v}{\xi\omega-\overline{l}}\right)\xi(x,t)^k dx\\
\nonumber
&\le &\gamma\frac{(\xi\omega)^{p-2}}{r_0^{n+p}}\iint\limits_{L_0(\overline{l})}\left(1+\frac{\xi\omega-v}{\xi\omega-\overline{l}}\right)
^{(1-\lambda)(p-1)}\left(\frac{\xi\omega-v}{\xi\omega-\overline{l}}\right)^{2\lambda(p-1)}
\xi(x,t)^{k-p} dxdt, \\
\nonumber
 &+&\gamma\varepsilon\frac{(\xi\omega)^{p-2}}{r_0^{n+p}}\iint\limits_{L_0(\overline{l})}
 \left(1+\frac{\xi\omega-v}{\xi\omega-\overline{l}}\right)^{(1+\lambda)(p-1)}\xi(x,t)^{k-p}dxdt\\
\nonumber
&+&\gamma {\varepsilon}^{-\frac
1{p-1}}\frac{r_0^{p-n}}{(\xi\omega-\overline{l})^p}\int
\limits_{B_0}F_1(x) dx +\gamma
\frac{r_0^{p-n}}{(\xi\omega-\overline{l})^{p-1}}\int\limits_{B_0}
F_2(x) dx\\
&\le& \gamma\nu \theta(\xi\omega)^{p-2}C^{n+p}+\gamma
{\varepsilon}^{-\frac
1{p-1}}C^{n-p}(B^{1-p}+B^{-p}).\label{e4.10}
\end{eqnarray}
Fix a number $\varkappa\in (0,1)$ depending on the known data.
First, choose $\varepsilon=\nu$, next choose $\nu$ from  the condition
$\gamma\nu\theta(\xi\omega)^{p-2}C^{n+p}\le \frac {\varkappa}8$
and $B$ from the condition $B^{1-p}\gamma\nu^{-\frac
1{p-1}}C^{n-p}\le \frac{\varkappa}4$. Then we obtain from \eqref{e4.9},
\eqref{e4.10} that $A_0(\overline{l})\le \frac{\varkappa}2$.

\begin{lemma}
\label{lem4.3}
Suppose we have chosen
$l_1, \dots, l_j$ and $\delta_0, \dots, \d_{j-1}$ such that
\begin{equation}\label{e4.11}
\frac{l_{i-1}}2+\frac B4 \alpha_{i-1}<l_i\le l_{i-1}-\frac
14(\alpha_{i-2}-\alpha_{i-1}), \quad i=1,2,\dots, j,
\end{equation}
\begin{equation}\label{e4.12}
A_{i-1}(l_i)\le \varkappa, \quad i=1,2,\dots, j,
\end{equation}
\begin{equation}\label{e4.13}
l_i>\frac {B{\alpha_{i-1}}}2, \quad i=1,2,\dots j.
\end{equation}
Then
\begin{equation}\label{e4.15}
A_j(\overline{l})\le \frac \varkappa 2, \quad \overline{l}=\frac
{l_j}2+\frac B4\alpha_{j}.
\end{equation}
\end{lemma}
\proof
Let us decompose $L_j(\overline{l})$ as
\[
L_j(\overline{l}) =L_j\pr(\overline{l}) \cup L_j^{\prime\prime}(\overline{l}), \quad L_j\pr(\overline{l})=\left\{\frac{l_j-v}{l_j-\overline{l}}
\le\varepsilon_1\right\}, \quad L_j^{\prime\prime}(\overline{l})=L_j(\overline{l})\setminus L_j\pr(\bar l).
\]
Using that $\xi_{j-1}(x,t)=1$ for $(x,t)\in
Q_j(\overline{l})$ and inequality \eqref{e4.12}
we have
\begin{eqnarray}
\nonumber
&&\frac{(l_j-\overline{l})^{p-2}}{r_j^{n+p}}\iint\limits_{L'_j(\overline{l})}\left(\frac{l_j-v}{l_j-\overline{l}}
\right)^{(1+\lambda)(p-1)}\xi_j(x,t)^{k-p}dx dt\\
\nonumber
&\le& \frac{(l_j-\overline{l})^{p-2}}{r_j^{n+p}}\varepsilon_1^{(1+\lambda))(p-1)}|L_j(\overline{l})|
\le \frac{\varepsilon_1^{(1+\lambda)(p-1)}}{r_j^n}
\esssup_t\int\limits_{L_j(\overline{l},t)}\xi_{j-1}(x,t)^k dx\\
&\le& \frac{\varepsilon_1^{(1+\lambda)(p-1)}}{r_j^n}
\esssup_t\int_{L_{j-1}(\bar l,t)} G\left(\frac{l_{j-1}-v}{l_{j-1}-l_j}\right)\xi_{j-1}(x,t)^kdx
\le 2^n \eps_1^{(1+\l)(p-1)}\varkappa.
\label{e4.16}
\end{eqnarray}
Above we also used the following inequality, which follows from \eqref{e4.11}, \eqref{e4.13}.
\begin{eqnarray}
\nonumber
l_j-\bar l=\frac{l_j}2-\frac{B}4\a_{j}\ge \frac{l_{j-1}}4+\frac{B}8\a_{j-1}-\frac{B}4\a_{j}\\
\label{e4.17}
=\frac{l_{j-1}-l_j}4+\frac{l_j}4+\frac{B}8\a_{j-1}-\frac{B}4\a_{j}\ge \frac{l_{j-1}-l_j}4.
\end{eqnarray}
It follows from \eqref{e4.17} that $Q_j(\bar l)\subset Q_{j-1}(l_j)$.

Using the Young inequality we have

\begin{eqnarray}
\nonumber
&&\frac{(l_j-\bar{l})^{p-2}}{r_j^{n+p}}\iint\limits_{L_j^{\prime\prime}(\bar{l})} \left(\frac{l_j-v}{l_j-\bar{l}}\right)^{(1+\lambda)(p-1)}
\xi_j(x,t)^{k-p}dx dt\\
&\le& \varepsilon_1\frac{(l_j-\bar l)^{p-2}}{r_j^{n+p}}|L_j(\bar l)|+
\gamma(\varepsilon_1)\frac{(l_j-\bar l)^{p-2}}{r_j^{n+p}}
\iint\limits_{L_j^{\prime\prime}(\bar l)}\left(\frac{l_j-v}{l_j-\bar
l}\right)^{p\frac{n+\rho(\lambda)}{\rho(\lambda)n}}\xi_j(x,t)^{(k-p)z}dx
dt,
\label{e4.18}
\end{eqnarray}
where $\rho(\lambda)=\frac p{p-1-\lambda}, \
z(1+\lambda)(p-1)=p\frac{n+\rho(\lambda)}{\rho(\lambda)n}$ due to our choice of
$\lambda,\, z>1$.

Similarly to \eqref{e4.16}, the first term in the right hand side of  \eqref{e4.18} is estimated as
\begin{equation}
   \label{4.19}
   \varepsilon_1 \frac {(l_j-\bar l)^{p-2}}{r_j^{n+p}}|L_j(\bar l)|\le \varepsilon_12^n\varkappa.
\end{equation}
Define
\[
w_j(x,t)=\frac 1{l_j-\bar l}\left[\int^{l_j}_{v(x,t)} \left(1+\frac {l_j-s}{l_j-\bar l}\right)^{-\frac 1p+\frac \lambda p}\left(\frac{l_j-s}{l_j-\bar l}\right)^{-\frac{2\lambda} p}ds\right]_+.
\]
Using the embedding theorem and Lemma~\ref{lem3.2} we have
\begin{eqnarray}
\nonumber
&&\gamma(\varepsilon_1)\frac {(l_j -\bar l)^{p-2}}{r_j^{n+p}}\iint\limits_{L_j^{\prime\prime}(\bar l)}
\left(\frac{l_j-v}{l_j-\bar l}\right)^{p\frac{n+p(\lambda)}{\rho(\lambda) n}}\xi_j(x,t)^{(k-p)z} dx dt\\
\nonumber
&\le& \gamma(\varepsilon_1)\frac{(l_j-\bar l)^{p-2}}{r_j^{n+p}} \iint\limits_{L_j^{\prime\prime}}
w_j(x,t)^{p\frac{n+p(\lambda)}{n}} \xi_j(x,t)^{(k-p)z}dx dt\\
\nonumber
&\le& \gamma(\varepsilon_1)\frac{(l_j -\bar l)^{p-2}}{r_j^{n+p}}
\left(\esssup \limits_t \iint\limits_{L_j^{\prime\prime}(\bar l, t)}w_j(x,t)^{\rho(\lambda)}
\xi_j (x,t)^{\frac{\rho(\lambda)(k-p)}{p(n+p(\lambda))}} dx\right)^{\frac pn}
\end{eqnarray}
\begin{eqnarray}
\nonumber
&\times& \iint\limits_{L_j(\bar l, t)}\left|\n \left(w_j(x,t)\xi_j^{\frac{\rho(\lambda)(k-p)zn}{p(n+p(\lambda))}}\right)\right|^p dxdt\\
\nonumber
&\le&\gamma(\varepsilon_1)\biggl\{\frac{(l_j-\bar l)^{p-2}}{r_j^{n+p}}\iint\limits_{L_j(\bar l)}\left(1+\frac{l_j-v}{l_j-\bar l}\right)^{(1-\lambda)(p-1)}\left(\frac{l_j-v}{l_j-\bar l}\right)^{2\lambda(p-1)}dx dt\\
\nonumber
&+& \gamma\varepsilon\frac{(l_j-\bar l)^{p-2}}{r_j^{n+p}}\iint\limits_{L_j(\bar l)}\left(1+\frac {l_j-v}{l_j-\bar l}\right)^{(1+\lambda)(p-1)}dxdt\\
&+& \gamma\varepsilon^{-\frac1{p-1}}\frac{r_j^{p-n}}{(l_j-\bar l)^p}\int_{B_j} F_1(x) dx
+\gamma\frac{r_j^{p-n}}{(l_j-\bar l)^{p-1}}\int_{B_j} F_2(x) dx \biggr\}^{1+\frac pn}.
\label{e4.20}
\end{eqnarray}
Let us take $\varepsilon=1$. Using the inequality $l_j\le l_{j-1}$ and \eqref{e4.10}, \eqref{e4.12}, \eqref{e4.17}
we have
\begin{eqnarray}
\nonumber
&&\frac{(l_j-\bar l)^{p-2}}{r_j^{n+p}}\iint\limits_{L_j(\bar l)}\left(1+\frac{l_j-v}{l_j-\bar l}\right)^{(1-\lambda)(p-1)}
\left(\frac{l_j-v}{l_j-\bar l}\right)^{2\lambda(p-1)}dx dt\\
\nonumber
&+& \frac{(l_j-\bar l)^{p-2}}{r_j^{n+p}}\iint\limits_{L_j(\bar l)}\left(1+\frac{l_j-v}{l_j-\bar l}\right)^{(1-\lambda)(p-1)} dx dt\\
\nonumber
&\le&\gamma\frac{(l_j-\bar l)^{p-2}}{r_j^{n+p}}|L_j(\bar l)| +\gamma\frac{(l_j-\bar l)^{-1-\lambda(p-1)}}{r_j^{n+p}}
\iint\limits_{L_j(\bar l)}(l_{j-1}-v)^{(1+\lambda)(p-1)}\xi_{j-1}(x,t)^{k-p} dx dt\\
&\le& \gamma \varkappa +\gamma \frac{\delta^{p-2}_{j-1}}{r_j^{n+p}}\iint\limits_{L_{j-1}(\bar l)}
\left(\frac{l_{j-1}-v}{l_{j-1}-l_j}\right)^{(1+\lambda)(p-1)} \xi_{j-1}(x,t)^{k-p} dx dt
\le \gamma \varkappa.
\label{e4.21}
\end{eqnarray}
Furthermore,  \eqref{e4.13} implies that
$l_j-\bar l\ge \frac B4 (\alpha_{j-1}-\alpha_j) 
$.
Therefore by \eqref{alpha-alpha} we have
\begin{equation}
\gamma\frac{r_j^{p-n}}{(l_j-\bar l)^p}\int_{B_j} F_1(x) dx +\gamma\frac{r_j^{p-n}}{(l_j-\bar l)^{p-1}}\int_{B_{j}}F_2(x) dx \le \gamma(B^{1-p}+B^{-p}).
\label{e4.22}
\end{equation}
Using Lemma \ref{lem3.2} again, we obtain
\begin{eqnarray}
\nonumber &&\frac 1{r_j^n}\int\limits_{L_j(\bar l, t)} G\left(\frac{l_j-v}{l_j-\bar l}\right) \xi_j (x,t)^k dx \\
\nonumber &\le&\gamma\frac{(l_j-\bar l)^{p-2}}{r_j^{n+p}}\iint\limits_{L_j(\bar l)}
\left(1+\frac{l_j-v}{l_j-\bar l}\right)^{(1-\lambda)(p-1)}\left(\frac{l_j-v}{l_j-\bar l}\right)^{2\lambda(p-1)}\xi_j(x,t)^{k-p} dx dt\\
\nonumber
&+& \gamma\varepsilon\frac{(l_j-\bar l)^{p-2}}{r_j^{n+p}}\iint\limits_{L_j(\bar l)}
\left(1+\frac{l_j-v}{l_j-\bar l}\right)^{(1+\lambda)(p-1)}\xi_j(x,t)^{k-p} dx dt\\
&+&\gamma \varepsilon^{-\frac 1{p-1}}\frac{r_j^{p-n}}{(l_j-\bar l)^p} \int_{B_j} F_1(x) dx
+\gamma\frac{r_j^{p-n}}{(l_j-\bar l)^{p-1}} \int_{B_j} F_2(x) dx.
\label{e4.23}
\end{eqnarray}
Using the decomposition $L_j(\bar l)=L^\prime (\bar l) \cup L^{\prime\prime}(\bar l)$
we have
\begin{eqnarray}
\nonumber
&&\frac{\gamma(l_j-\bar l)^{p-2}}{r_j^{n+p}}\iint\limits_{L_j(\bar l)}\left(1+\frac {l_j-v}{l_j-\bar l}\right)^{(1-\lambda)(p-1)}
\left(\frac{l_j-v}{l_j-\bar l}\right)^{2\lambda(p-1)}\xi_j(x,t)^{k-p}  dx dt\\
\nonumber
&\le& \frac{\gamma(l_j-\bar l)^{p-2}}{r_j^{n+p}}\varepsilon_1^{2\lambda(p-1)}(1+\varepsilon_1)^{(1-\lambda)(p-1)} |L_j^\prime(\bar l)|\\
&+& \gamma(\varepsilon_1) \frac{(l_j-\bar l)^{p-2}}{r_j^{n+p}}\iint\limits_{L_j^{\prime\prime}(\bar l)}
\left(\frac{l_j-v}{l_j-\bar l}\right)^{(1+\lambda)(p-1)} \xi_j(x,t)^{k-p}  dx dt.
\label{e4.24}
\end{eqnarray}
Similarly
\begin{eqnarray}
\nonumber
&&\gamma\varepsilon\frac{(l_j-\bar l)^{p-2}}{r_j^{n+p}}\iint\limits_{L_j(\bar l)}
\left(1+\frac{l_j-v}{l_j-\bar l}\right)^{(1-\lambda)(p-1)}\xi_j(x,t)^{k-p}  dx dt\\
\nonumber
&\le& \gamma\varepsilon(1+\varepsilon_1)^{(1-\lambda)(p-1)} \frac {(l_j-\bar l)^{p-2}}{r_j^{n+p}}|L^\prime_j(\bar l)|\\
&+& \varepsilon\gamma(\varepsilon_1)\frac{(l_j-\bar l)^{p-2}}{r_j^{n+p}}
\iint\limits_{L_j^{\prime\prime}(\bar l)}\left(\frac{l_j-v}{l_j-\bar l}\right)^{(1+\lambda)(p-1)}\xi_j(x,t)^{k-p} dx dt.
\label{e4.25}
\end{eqnarray}

Combining estimates \eqref{e4.16}-\eqref{e4.25} we have
\begin{equation}
\label{e4.26}
A_j(\bar l) \le
\gamma(\varepsilon_1^{(1+\lambda)(p-1)}+\varepsilon_1^{2\lambda(p-1)}+\varepsilon(1+\varepsilon_1)^{(1-\lambda)(p-1)}) \varkappa
+ \gamma(\varepsilon_1,\eps)(B^{1-p}+B^{-p})
+ \gamma(\varepsilon_1,\eps)\{\varkappa+(B^{1-p}+B^{-p})\}^{1+\frac pn.}
\end{equation}
First choose  $\varepsilon_1$ from the condition
\begin{equation}\label{e4.27}
\varepsilon_1^{(1+\lambda)(p-1)}+\varepsilon_1^{2\lambda(p-1)} =\frac 1{16}.
\end{equation}
Next we choose $\varepsilon$ from the equality
\begin{equation}\label{e4.28}
    \varepsilon(1+\varepsilon_1)^{(1-\lambda)(p-1)}=\frac 1{16}.
\end{equation}
Fix $\varkappa$ by
\begin{equation}\label{e4.30}
    \gamma(\varepsilon_1, \varepsilon)\varkappa^{\frac pn}=\frac 1{16},
\end{equation}
and choosing $B$ large enough so that
\begin{equation}\label{e4.29}
    B^{1-p}+B^{-p}\le \frac{\varkappa}{16},
\end{equation}
we conclude
from \eqref{e4.26} that
$A_j(\bar l)\le \frac \varkappa 2$, which completes the proof of Lemma \ref{lem4.3}.\qed

\bigskip
Further, since $A_j(l)$ is an increasing and continuous function and $A_j(l)\to \infty$ if $l\to l_j$,
inequality \eqref{e4.15} ensures the existence of $\tilde{l} \in (\bar l, l_j)$ such that  $A_j(\tilde{l})=\varkappa$.
If $\tilde l< l_j-\frac {1}4(\alpha_{j-1}-\alpha_{j})$ we set $l_{j+1}=\tilde l$.
If  $\tilde l\ge l_j-\frac{1}4(\alpha_{j-1}-\alpha_{j})$, then we set
$l_{j+1}=l_j-\frac {1}4(\alpha_{j-1}-\alpha_{j})$ and in both cases we set $\delta_j=l_j-l_{j+1}$.

In what follows
\[Q_j=Q_j(l_{j+1}), \quad L_j=L_j(l_{j+1})\]
\begin{lemma}\label{lem4.4}
Let the conditions of Theorem~\ref{thm1.4} be fulfilled. Then for any $j\ge 1$ the following inequality holds
\begin{equation}\label{e4.31}
\delta_j\le \frac 12 \delta_{j-1}+
\gamma r_j
+\gamma\left(r_j^{p-n}\int_{B_{j-1}} F_1(x) dx\right)^{\frac 1p}+\gamma\left(r_j^{p-n}\int_{B_{j-1}}F_2(x) dx\right)^{\frac 1{p-1}}.
\end{equation}
\end{lemma}
\proof
Fix $j\ge 1$ and assume without loss that
\begin{equation}\label{e4.32}
\delta_j> \frac 12 \delta_{j-1}, \quad \delta_j> \frac{1}4(\alpha_{j-1}-\alpha_j)
\end{equation}
since in the opposite case
due to \eqref{alpha-alpha1}
inequality \eqref{e4.31}  is obvious.  The second inequality in \eqref{e4.32}
ensures that  $A_j(l_{j+1})=\varkappa$. Using the decomposition
$L_j=L_j^\prime\cup L_j^{\prime\prime}$  similarly to  \eqref{e4.16}, \eqref{e4.18}--\eqref{e4.21} we obtain
\begin{eqnarray}
\nonumber
&&\frac{\delta_j^{p-2}}{r_j^{n+p}}\iint\limits_{L_j} \left(\frac{l_j-v}{\delta_j}\right)^{(p-1)(1+\lambda)}\xi_j(x,t)^{k-p} dx dt\\
\nonumber
&\le&\varepsilon_1^{(1+\lambda)(p-1)}\frac{\delta_j^{p-2}}{r_j^{n+p}}|L_j^\prime|+\frac{\delta_j^{p-2}}{r_j^{n+p}}\iint\limits_{L_j^{\prime\prime}}
\left(\frac{l_j-v}{\delta_j}\right)^{(1+\lambda)(p-1)}\xi_j(x,t)^{k-p} dx dt\\
&\le&\gamma(\varepsilon_1^{(1+\lambda)(p-1)}+\varepsilon_1)\varkappa +\gamma(\varepsilon_1)
\biggl\{\varkappa+\delta^{-p}_jr_j^{p-n}\int\limits_{B_j}F_1(x) dx
+\delta_j^{1-p} r_j^{p-n}\int_{B_j} F_2(x) dx\biggr\}^{1+\frac pn},
\label{e4.33}
\end{eqnarray}
Using Lemma \ref{lem3.2}, in the same way as
\eqref{e4.23} -\eqref{e4.25} we have

\begin{eqnarray}
\nonumber
&&\esssup\limits_t\frac 1{r_j^n}\int_{L_j(t)}G\left(\frac{l_j-v}{\delta_j}\right)\xi_j(x,t)^k dx
\nonumber \le\gamma(\varepsilon_1^{(1+\lambda)(p-1)}+\varepsilon_1^{\lambda(p-1)}+\varepsilon(1+\varepsilon_1)^{\lambda(p-1)})\varkappa\\
\nonumber
&+& \gamma(\varepsilon_1, \varepsilon)\left[\delta_j^{-p}r_j^{p-n}\int\limits_{B_j} F_1(x) dx+\delta_j^{1-p}r_j^{p-n}\int\limits_{B_j}F_2(x) dx\right]\\
&+& \gamma(\varepsilon_1, \varepsilon)\left\{ \varkappa+\delta_j^{-p}r_j^{p-n}\int\limits_{B_j} F_1(x) dx+\delta_j^{1-p}r_j^{p-n}\int\limits_{B_j}F_2(x) dx\right\}^{1+\frac pn}.
\label{e4.34}
\end{eqnarray}

Choosing $\varepsilon_1, \varepsilon, \varkappa$ from inequalities \eqref{e4.27}, \eqref{e4.28}, \eqref{e4.30} we conclude that
at least one of  the two following inequalities holds
\[\delta_j\le \gamma\left(r_j^{p-n}\int_{B_j} F_1(x) dx\right)^{\frac 1p},\quad
\delta_j\le \gamma\left(r_j^{p-n}\int\limits_{B_j} F_2(x) dx\right)^{\frac 1{p-1}},\]
which proves Lemma \ref{lem4.4}. \qed

Summing up  inequality \eqref{e4.31} with respect to $j=1, \dots, J-1$ we obtain
\begin{equation}l_1-l_J\le \delta_0 +
\gamma r_0
+\gamma\sum^\infty_{j=1}\left(r_j^{p-n}\int\limits_{B_{j-1}} F_1(x) dx\right)^{\frac 1p}+\gamma\sum^\infty_{j=1}\left(r_j^{p-n}\int\limits_{B_{j-1}} F_2(x) dx\right)^{\frac 1{p-1}}.
\label{e4.35}
\end{equation}
If $l_1$ is defined  by $l_1=\xi\omega-\frac 14(\alpha_{-1}-\alpha_0)$ then $\delta_0=\frac 14(\alpha_{-1}-\alpha_0)$.
Passing  to the limit in \eqref{e4.35} as $J\to \infty$ we have
\begin{equation}
\xi\omega\le \lim_{j\to \infty} l_j
+\gamma r_0
+\gamma\int_0^\rho\left(\frac 1{r^{n-p}}\int\limits_{B_r(y)}F_1(x)dx\right)^{\frac 1p} \frac{dr}r
+ \gamma\int\limits^\rho_0\left(\frac 1{r^{n-p}}\int\limits_{B_{r}(y)}F_2(x)dx\right)^{\frac 1{p-1}} \frac{dr}r.
\label{e4.36}
\end{equation}
If $l_1<\xi\omega-\frac 14(\alpha_{-1}-\alpha_0)$ then $A_0(l_1)=\varkappa$ and  at least one of the following inequalities holds
\begin{equation}\label{e4.37}
\frac{\delta_0^{p-2}}{r_0^{n+p}}\iint\limits_{L_{0}}\left(\frac{l_0-v}{\delta_0}\right)^{(1+\lambda)(p-1)} dx dt \ge \frac \varkappa 2
\end{equation}
or
\begin{equation}\label{e4.38}
   \esssup\limits_t \frac 1{r_0^n} \int\limits_{L_{0}(t)}G\left(\frac{l_0-v}{\delta_0}\right) dx\ge \frac \varkappa 2.
\end{equation}
Similarly to \eqref{e4.9}, it follows from \eqref{e4.37} that
\begin{equation}\label{e4.39}
\frac\varkappa 2\le \frac{\delta_0^{p-2}}{r_0^{n+p}}\iint\limits_{L_0}\left(\frac{l_0-v}{\delta_0}\right)^{(1+\lambda)(p-1)} dx dt\le \gamma C^{n+p}\frac{\nu(\xi\omega)^{(1+\lambda)(p-1)}\theta}{\delta_0^{1+\lambda(p-1)}}.\end{equation}
Similarly to \eqref{e4.10} with $\varepsilon=\nu^{-\frac12}$, it follows from \eqref{e4.38} that
\begin{eqnarray}
\nonumber
&&\frac \varkappa 2\le \esssup\limits_t \frac 1{r_0^n} \int\limits_{L_0(t)}G\left(\frac{l_0-v}{\delta_0}\right) dx
\\
\nonumber
&&
\le
\gamma\frac{\delta_0^{p-2}}{r_0^{n+p}} \nu^{-\frac12}\iint\limits_{L_0}\left(1+\frac{l_0-v}{\delta_0}\right)^{(1+\lambda)(p-1)} dx dt
+\gamma\frac{\nu^{\frac1{2(p-1)}}r_0^{p-n}}{\delta_0^p}\int\limits_{B_0} F_1(x) dx+\gamma\frac{r_0^{p-n}}{\delta_0^{p-1}}\int\limits_{B_0} F_2(x) dx\\
&&\le
\gamma\frac{C^{n+p}\nu^\frac12(\xi\omega)^{(1+\lambda)(p-1)}\theta}{\delta_0^{1+\lambda(p-1)}}
+\gamma\frac{\nu^{\frac 1{2(p-1)}}}{\delta_0^p} r_0^{p-n}\int\limits_{B_0} F_1(x) dx +\gamma\frac{r_0^{p-n}}{\delta_0^{p-1}}\int\limits_{B_0} F_2(x) dx.
\label{e4.40}
\end{eqnarray}

First we choose $C>16$.
Then \eqref{e4.39}, \eqref{e4.40} imply that
\begin{eqnarray}
\nonumber
&&\delta_0\le \xi\omega(\gamma\varkappa^{-1}C^{n+p}\nu(\xi\omega)^{p-2}\theta)^{\frac 1{1+\lambda(p-1)}}\\
&+&\gamma\nu^{-\frac1{p(p-1)}}\left(r_0^{p-n}\int\limits_{B_0} F_1(x) dx\right)^{\frac 1{p}}
+\gamma\left(r_0^{p-n}\int\limits_{B_0}F_2(x) dx\right)^{\frac1{p-1}}.
\label{e4.41}
\end{eqnarray}
Finally due to the inequality $\xi\omega\ge \theta^{-\frac 1{p-2}}$ from \eqref{e4.35}, \eqref{e4.36}, \eqref{e4.41} we  have due to $\nu<1$
\begin{eqnarray}
\nonumber
&&\xi\omega\le v(x_1, t_1)+ \xi\omega(\gamma \varkappa^{-1}C^{n+p}\nu^\frac12 (\xi\omega)^{p-2} \theta)^{\frac 1{1+\lambda(p-1)}}+\g \rho \\
&+&
\g \int\limits_0^\rho\left(r^{p-n}
\int\limits_{B_r(y)}F_1(x) dx\right)^{\frac 1p} \frac{dr}r
+
\g \int\limits_0^\rho \left(r^{p-n}\int\limits_{B_r(y)} F_2(x) dx\right)^{\frac 1{p-1}}\frac {dr}r.
\label{e4.42}
\end{eqnarray}

Next we fix $\nu$ from the condition
\begin{equation}
\gamma\varkappa^{-1}C^{n+p}\nu(\xi\omega)^{p-2} \theta=\left(\frac{1-a}2\right)^{1+\lambda(p-1)}
\label{e4.44}
\end{equation}
and finally, choosing $B$ large enough so that
\begin{equation}\label{e4.45}
B\ge
\frac{2\g}{1-a},
\end{equation}
we obtain from \eqref{e4.42}
\begin{equation}\label{e4.46}
u(x_1,t_1)\ge \mu_-+a\xi\omega.
\end{equation}
Since $(x_1, t_1)$ is an arbitrary point in $Q_{\rho}^\theta(y,s)$,
from \eqref{e4.46} the required \eqref{e1.13} follows, which proves
Theorem \ref{thm1.4}.\qed

 \section{Expansion of positivity. Proof of Theorem~\ref{thm1.5}}
 \label{expansion}

In the proof we closely follow \cite{DiGV}, also using the idea of logarithmic estimates from \cite{DiB}.
Our assumption
here are again \eqref{e1.2}.
 In what follows we suppose that
 \begin{equation}\label{e5.1}
 N\ge B(\rho+\calF_1(2\rho)+\calF_2(2\rho)).
 \end{equation}

 Let $0\le \tau\le \frac 12(p-2) \ln B$, $k=\mu_-+e^{-\frac\tau{p-2}} N, \ \theta=e^\tau N^{-(p-2)}$,
 $\xi\in C_0^\infty(B_\rho(y)), \xi(x)= 1$ if $x\in B_{\rho/2}(y)$,
 $0\le\xi(x)\le 1,  \left|{\n\xi(x)}\right|\le 2\rho^{-1}$.
As above, set
 $\Psi_-(u)=\ln_+\frac{H^-_k}{H^-_k-(k-u)_++2^{-s_0}e^{-\frac{\tau}{p-2}}N}, \ $
 $s_0$
 is a positive number satisfying $s_0< \frac 12 \ln B$, which will be determined later
 depending on the data.
 Note the evident inequalities
 \[
 (k-u)_+\le e^{-\frac{\tau}{p-2}}N, \quad |\Psi_-(u)|\le s_0\ln 2,
 \]
 \[
 |\Psi^\prime_-(u)|\le 2^{s_0} e^\frac{\tau}{p-2} N^{-1},\quad |\Psi^\prime_-(u)|^{2-p}\le e^{-\tau} N^{p-2},
 \]
 \[
 \int\limits_{B_\rho(y)}F_1(x) dx\le \gamma \calF_1(2\rho)^p\rho^{n-p}, \int\limits_{B_\rho(y)} F_2(x) dx\le\gamma \calF_2(2\rho)^{p-1}\rho^{n-p}.
 \]

 Since
  condition \eqref{e1.16} guarantees that
  \begin{equation}\label{e5.2}
  \Psi_-(u)=0\   \mbox{for}\ x\in B_\rho(y), t=s,
  \end{equation}
  Lemma \ref{lem2.3} implies that
  \begin{eqnarray}
  \nonumber
  &&\esssup\limits_{s<t<s+\theta\rho^p}\int\limits_{B_{\frac \rho 2}(y)}\Psi^2_-(u)dx
  \le \gamma\iint\limits_{Q^\theta_\rho(y,s)}\Psi_-|\Psi_-^\prime(u)|^{2-p}|\n\xi|^p dx dt
  +\gamma \iint\limits_{Q^\theta_\rho(y,s)}\Psi_-|\Psi_-^\prime(u)|^{2-p} dx dt\\
  \nonumber
  &+&\gamma\iint\limits_{Q^\theta_\rho(y,s)}\Psi_-|\Psi_-^\prime(u)|^2 F_1(x) dx dt
  +\iint\limits_{Q^\theta_\rho(y,s)}\Psi_-|\Psi_-^\prime(u)|F_2(x) dx dt\\
  &\le& \gamma s_0\rho^n
  +
  \gamma s_0\left(\frac{2^{s_0}e^\frac{\tau}{p-2} \calF_1(2\rho)}N\right)^p\rho^n+\gamma s_0\left(\frac {2^{s_0}e^\frac{\tau}{p-2} \calF_2(2\rho)} N \right)^{p-1} \rho^n
  \le \gamma s_0 \rho^n.
  \label{e5.3}
\end{eqnarray}
Since
\[\Psi_-(u) \ge (s_0-1)\ln 2\  \mbox{for}\ x\in B_{\rho/2}(y)\cap\left\{u<\mu_-+\frac{e^{-\frac\tau{p-2}}N}{2^{s_0}}\right\}, \]
 inequality \eqref{e5.3} yields
\begin{equation}
\label{e5.4}
\left|\left\{ x\in B_{\rho/2}(y): u{(x,t)}<\mu_-+\frac{e^{-\frac\tau{p-2}}N}{2^{s_0}}\right\}\right|\le\gamma\frac {s_0}{{(s_0-1)^2}}| B_{\rho/2}(y) |
\end{equation}
for all $t\in (s,s+\theta \rho^p), \ 0\le \tau\le \frac12 (p-2) \ln B$.
Choosing $s_0$ from the condition
\begin{equation}\label{e5.5}
\gamma\frac {s_0}{(s_0-1)^2}\le \frac 12,
\end{equation}
we obtain
\begin{equation}
\label{e5.6}
\left|\left\{ x\in
B_{\rho/2}(y): u(x,s+e^\tau N^{-(p-2)}\rho^p)\le
\mu_-+\frac{e^{-\frac\tau{p-2}}N}{2^{s_0}}\right\}\right|\le \frac 12
|B_{\rho/2}(y)|
\end{equation}
for all $0\le \tau\le \frac 12 (p-2)\ln B$.

In the same way as in \cite[p.\,191]{DiGV} we consider the function
\[
w(x,\tau)=e^\frac\tau{p-2} N^{-1} \rho^{\frac p{p-2}}\left(u(x,s+(e^\frac\tau{p-2} N^{-1})^{p-2}\rho^p)-\mu_-\right),
\]
and let $k_0=2^{-s_0} \rho^{\frac p{p-2}}$.

Inequality \eqref{e5.6} translates into $w$ as $|\{x\in B_{\rho/2}(y) :  w(x,\tau)\le k_0\}|\le \frac12 |B_{\rho/2}(y)|$, which
yields
\begin{equation}
\label{e5.7}
|\{x\in B_{4\rho} (y): w(x,\tau)\le k_0\}|\le \left(1-\frac 1{2\cdot 8^n}\right)|B_{4\rho}(y)|
\end{equation}
for all $\tau\in (0, 2^{-s_0(p-2)-1}(p-2)\ln B\, k_0^{2-p}\rho^p)$.

Since $w\ge 0$, formal differentiation, which can be justified in a standard way, gives
\begin{equation}
\label{e5.8} w_\tau=\frac1{p-2}w+(e^\frac\tau{p-2} N^{-1}\rho^{\frac p{p-2}})^{p-1}u_t\ge {\rm div}\,\tilde A(x,t,w,\n w)+\tilde{a}_0(x,t,w,\n w),
\end{equation}
where $\tilde{A}, \ \tilde a_0$ satisfy the inequalities

\begin{eqnarray}
\nonumber
\tilde{A}(x,t,w,\n w)\cdot\n w&\ge& c_1|\n w|^p,\\
\nonumber
|\tilde{A}(x,t,w,\n w)|&\le& c_2|\n w|^{p-1}+(e^\frac\tau{p-2} N^{-1}\rho^{\frac p{p-2}})^{p-1}f_1(x),\\
\label{e5.9}
|\tilde{a}_0(x,t,w,\n w)|&\le& c_2 |\n w|^{p-1}+(e^\frac\tau{p-2} N^{-1} \rho^{\frac p{p-2}})^{p-1}f_2(x).
\end{eqnarray}

\begin{lemma}\label{lem5.1}
For every $\nu\in (0,1)$ there exists $s_*>s_0$\,,\ $2^{s_\ast}\le 2^{-s_0( p-2)-1}(p-2)\ln B$, depending only on the data and $\nu$ such that
\begin{equation}\label{e5.10}
|\{Q^\ast_{4\rho}:w(x,\tau)<\frac{k_0}{2^{s_\ast}}\}|\le \nu|Q^\ast_{4\rho}|,
\end{equation}
where $Q^\ast_\rho=B_\rho(y)\times(2k_0^{2-p}\rho^p, (2^{s_\ast}k_0^{-1})^{p-2}\rho^p).$
\end{lemma}
\proof
Using Lemma \ref{lem2.3m} with $k=\frac{k_0}{2^s},\, l=\frac{k_0}{2^{s-1}}, s_0\le s\le s_\ast$, due to \eqref{e5.7}
we obtain the inequality
\begin{equation}\label{e5.11}
\frac{k_0}{2^s}|A_{\frac{k_0}{2^{s}}, 4\rho}(\tau)|\le \gamma\rho\int\limits_{A_{\frac{k_0}{2^{s-1}}, 4\rho}(\tau)\setminus A_{\frac{k_0}{2^{s}},4\rho}(\tau)}\left|{\n w (x,\tau)}\right| dx
\end{equation}
for all $\tau\in(0, 2^{-s_0(p-2)-1}(p-2)\ln B k_0^{2-p}\rho^p)$, where
$A_{k,\rho}(\tau)=\{x\in B_\rho(y): w(x,\tau)\le k\}$.

Integrating the last inequality with respect to $\tau$, $\tau\in (2
k_0^{2-p}\rho^p, (2^{s_\ast}k_0^{-1})^{p-2}\rho^p)$, and using the
H\"older inequality we obtain
\begin{equation}
\label{e5.12}
\left(\frac{k_0}{2^s}\right)^{\frac p{p-1}}\left|A_{\frac{k_0}{2^{s}}, 4\rho}\right|^\frac{p}{p-1}\le \gamma
\rho^{\frac p{p-1}}
\left(\iint\limits_{A_{\frac{k_0}{2^{s-1}}, 4\rho}}\left|{\n w(x,\tau)}\right|^p dx dt\right)^{\frac 1{p-1}}
\left|A_{\frac{k_0}{2^{s-1}}, 4\rho}\setminus A_{\frac{k_0}{2^{s}}, 4\rho}\right|,
\end{equation}
where $A_{k,\rho}=\displaystyle\int^{(2^{s_\ast}k_0^{-1})^{p-2}\rho^p}_{2k_0^{2-p}\rho^p}A_{k,\rho}(\tau) d\tau.$

To estimate the first factor we use Lemma~\ref{lem2.2} with
$\displaystyle k=\frac{k_0}{2^s}$ and $\xi\in C^\infty_0 (\tilde{Q}^\ast_{8\rho}),\quad \xi(x,\tau)= 1$
for  $(x,t)\in Q^\ast_{4\rho},\quad 0\le \xi(x,\tau)\le 1,\quad \left|{\n \xi} \right|\le \gamma\rho^{-1},
\left|\frac{\partial\xi}{\partial \tau}\right|\le \gamma\left(\frac{k_0}{2^{s_\ast}}\right)^{p-2}\rho^{-p}$,\quad
$\tilde Q^\ast_\rho=B_\rho(y)\times \left( \frac{k_0^{2-p}}{4^{p+1}}\rho^p, (2^{s_\ast}k_0^{-1})^{p-2}\rho^p\right)$.

Due to \eqref{e5.8}, \eqref{e5.9} we obtain
\begin{eqnarray}
\nonumber
&&\iint\limits_{A_{\frac{k_0}{2^{s-1}},4\rho}}\left|{\n w(x,\tau)}\right|^pdx d\tau
\le \gamma\iint\limits_{\tilde Q^\ast_{8\rho}}\left[\left(\frac{k_0}{2^{s-1}}-w\right)^p_+\left|{\n\xi}\right|^p
+\left(\frac{k_0}{2^{s-1}}-w\right)^2_+\left|\frac{\partial \xi}{\partial t}\right|\right]dxd\tau\\
&&\ \ +
\gamma (e^\frac\tau{p-2} N^{-1}\rho^{\frac p{p-2}})^{p-1}\iint\limits_{\tilde Q^\ast_{8\rho}}(\frac{k_0}{2^{s-1}}-w)_+F_2(x) dx d\tau
+ \gamma(e^\frac\tau{p-2} N^{-1}\rho^{\frac p{p-2}})^p\iint\limits_{\tilde Q^\ast_{8\rho}}F_1(x) dx d\tau=\sum^3_{i=1} I_i.
\label{e5.13}\end{eqnarray}
Due to the choice of $\xi(x,t)$ we have
\begin{equation}\label{e5.14}
I_1\le \gamma\left(\frac{k_0}{2^s}\right)^p\rho^{-p} \left|Q^\ast_{4\rho}\right|.
\end{equation}
Using \eqref{e5.1} and the definition of the $K_p,\, \widetilde K_p$- classes we derive
\begin{equation}
I_2\le \gamma\left(e^\frac\tau{p-2} N^{-1}\rho^{\frac p{p-2}}\right)^{p-1}\frac{k_0}{2^s}\calF_2(2\rho)^{p-1} \rho^{-p}|Q^\ast_{4\rho}|=
\gamma\left(\frac {e^\frac\tau{p-2}2^{s_\ast}\calF_2(2\rho)}N\right)^{p-1}\left(\frac{k_0}{2^s}\right)^p\rho^{-p}|Q^\ast_{4\rho}|\le \gamma\left(\frac{k_0}{2^s}\right)^p\rho^{-p}|Q^\ast_{4\rho}|.
\label{e5.16}\end{equation}
\begin{equation}
I_3\le \gamma(e^\frac\tau{p-2} N^{-1}\rho^{\frac p{p-2}})^p \calF_1(2\rho)^p\rho^{-p}|Q^\ast_{4\rho}|
=\gamma\left(\frac{e^\frac\tau{p-2} 2^{s_\ast}\calF_1(2\rho)}{N}\right)^p\left(\frac{k_0}{2^s}\right)^p\rho^{-p}|Q^\ast_{4\rho}|
\le \gamma\left(\frac{k_0}{2^s}\right)^p\rho^{-p}|Q^\ast_{4\rho}|.
\label{e5.17}\end{equation}
Combining estimates \eqref{e5.12}--\eqref{e5.17} we obtain
\begin{equation}
\left|A_{\frac{k_0}{2^{s_*}},4\rho}
\right|^{\frac p{p-1}}\le \gamma|Q^\ast_{4\rho}|^{\frac 1{p-1}}\left|A_{\frac{k_0}{2^{s-1}}, 4\rho}\setminus A_{\frac{k_0}{2^{s}},4\rho}\right|.
\label{e5.18}
\end{equation}
Summing up the  last inequalities in $s$, $s_0<s\le s_\ast$, we conclude that
\begin{equation}\label{e5.19}
(s_\ast-s_0)|A_{\frac{k_o}{2^{s_\ast}},4\rho}|^{\frac p{p-1}}\le \gamma|Q^\ast_{4\rho}|^{\frac p{p-1}}.
\end{equation}
Choosing $s_\ast$ by the condition
\begin{equation}\label{e5.20}
(s_\ast-s_0)^{-\frac{p-1}p}\gamma\le \nu,
\end{equation}
we obtain inequality \eqref{e5.10}, which proves Lemma \ref{lem5.1}.\qed

\medskip
Using Theorem \ref{thm1.4} with $\xi=\displaystyle\frac 1{2^{s_\ast}}, \omega=k_0, \ \theta=(2^{s_\ast} k_0^{-1})^{p-2},\  a=\frac 12$ and choosing $\nu$
from condition~\eqref{e4.44}
we obtain
\begin{equation}
\label{e5.21}
w(x,\tau)\ge \frac N{2^{s_\ast+1}}\quad \mbox{for}\  \  x\in B_{2\rho}(y)
\end{equation}
and for all
$\tau\in \left(k_0^{2-p}(2\rho)^p, (2^{s_\ast}k_0^{-1})^{p-2} (2\rho)^p\right).$

Due to the choice of $k_0$, we have $ k_0^{2-p}\rho^p=2^{s_0(p-2)}$. For $\tau \in (k_0^{2-p}(2\rho)^p, (2^{s_\ast}k_0^{-1})^{p-2}(2\rho)^p)$
there holds
\[
\tilde b_1=\exp 2^{(s_0+1)(p-2)+1}\le e^\tau\le\exp 2^{(s_\ast+s_0+1)(p-2)+2}=\tilde b_2.
\]
Inequality \eqref{e5.21} translates for $u$ into
\begin{equation}
u(x,s+t)\ge 2^{-s_\ast-1}b_2^{-1} N=\sigma N\quad \mbox{ for} \ \ x\in B_{2\rho}(y)
\label{e5.22}\end{equation}
and for all $b_1N^{2-p}\rho^p\le t\le b_2N^{2-p}\rho^p,$ where
\begin{equation}
b_1=b_1(s_0)=\tilde b_1^{p-2}, b_2=b_2(s_0, s_\ast)=\tilde b_2^{p-2}
\label{e5.23}
\end{equation}
depend only on the data.
This completes the proof of Theorem \ref{thm1.5}.\qed

\section{Continuity of solutions. Proof of Theorem~\ref{thm1.1}}
\label{cont}

Here we closely follow \cite[Chapter\,III]{DiB}.
Let $(x_0, t_0)\in\Omega_T$ be arbitrary,
\[
Q_R(x_0, t_0)=B_R(x_0)\times(t_0-R^2, t_0), \quad R<\frac 14 \min \{1, t_0^{1/2}, {\rm dist}(x_0, \partial\Omega)\}.
\]
Set
\[\mu_+=\esssup\limits_{Q_R(x_0, t_0)}u(x,t),\qquad \mu_-=\essinf\limits_{Q_R(x_0, t_0)} u(x,t),\qquad \omega=\mu_+-\mu_-.\]
Fix a positive number $s^\ast$, $s_1=\frac{1}{p-2}\log_2b_1<s^\ast<\log_2b_2$, which will be determined later depending
 only on the known data, $b_1=b_1(s_0), b_2=b_2(s_0, s_\ast)$ are defined in \eqref{e5.23}.\\
 If
 \begin{equation}
 \omega \ge b_2(R+\calF_1(2R)+\calF_2(2R)),
 \label{e6.1}
 \end{equation}
 then the cylinder
 \[
 Q_R^\theta(x_0, t_0)=B_R(x_0)\times(t_0-\theta R^p, t_0), \, \theta=\left(\frac{2^{s^\ast}}\omega\right)^{p-2}\]
 is contained in $Q_R(x_0, t_0)$. In $Q_R^\theta(x_0, t_0)$ consider the cylinders
 \[
 Q^\eta_R(x_0,\bar t)=B_R(x_0)\times(\bar t-\eta R^p, \bar t),\qquad  \eta=b_1\omega^{2-p},\qquad t_0-\theta R^p\le \bar t-\eta R^p<\bar t\le t_0.
 \]
 Let us fix  $\nu\in (0,1)$  satisfying \eqref{e4.44} with $a=(\frac 12)^{\frac 1{p-2}}, \  \xi=\frac 12$ and $\theta=b_1\omega^{p-2}$.

The following two alternative cases are possible.

{\it First alternative.} There exists a cylinder $Q_R^\eta(x_0, \bar t)\subset Q^\theta_R(x_0, t_0)$ such that
\begin{equation}
\left|\left\{(x,t)\in Q^\eta_R(x_0, \bar t):u(x,t)\le \mu_-+\frac \omega2\right\}\right|\le \nu|Q^\eta_R(x_0, \bar t)|.
\label{e6.2}
\end{equation}

{\it Second alternative.} For all cylinders $Q_R^\eta(x_0, \bar t)\subset Q^\theta_R(x_0, t_0)$ the opposite inequality
\begin{equation}
 \left|\left\{(x,t)\in Q^\eta_R(x_0, \bar t):u(x,t)\le \mu_-+\frac \omega2\right\}\right|>\nu|Q_R^\eta(x_0, \bar t)|.
 \label{e6.3}
 \end{equation}
holds.

\subsection{Analysis of the first alternative}

Here we assume that \eqref{e6.1} is satisfied.
 By Theorem~\ref{thm1.4} with $\xi=\frac 12, a=\left(\frac 12\right)^{\frac 1{p-2}}$
 we obtain from \eqref{e6.2}\begin{equation}
 u(x,\bar t)\ge \mu_- +2^{-1-\frac 1{p-2}}\omega \quad \mbox{for all}\   x\in B_{\frac R4}(x_0).
 \label{e6.4}
 \end{equation}
 Using Theorem \ref{thm1.5} with $N=2^{-1-\frac 1{p-2}}\omega$  from \eqref{e6.4} we conclude that
 \begin{equation}
 u(x,t)\ge \mu_-+\sigma N\ \mbox{for}\  x\in B_{\frac R2}(x_0)
 \label{e6.5}\end{equation}
 and for all $t\in(\bar t+\frac 12 b_1 R^p, \bar t+\frac 12 b_2R^p)$,
 where $\sigma=\sigma(s_0, s_\ast),\, b_1=b_1(s_0),\, b_2=b_2(s_0, s_\ast)$ are fixed numbers defined in \eqref{e5.23}.

Since $\bar t+\frac 12 b_1R^p<t_0<\bar t+\frac 12b_2R^p$,  inequality~\eqref{e6.5} holds  for
 $(x,t)\in B_{\frac R2}(x_0)\times(t_0-\frac 12b_1R^p, t_0).$

 Thus we have proved the following
 \begin{proposition}\label{prop6.1}
 Suppose the first alternative holds. Then either
 $\omega\le b_2(R+\calF_1(2R)+\calF_2(2R))$, or \begin{equation}
 \essosc\limits_{Q^\eta_{\frac R2}(x_0, t_0)}u\le (1-\sigma)\omega.
 \label{e6.6}\end{equation}
 \end{proposition}
 \subsection{Analysis of the second alternative}

 This part is almost a literal repetition of the corresponding part from \cite[Chapter~III]{DiB} and is here for the readers' convenience.

 Since \eqref{e6.3} holds for all cylinders $Q^\eta_R(x_0, \bar t)$,  for the cylinders
  $Q_R^\eta(x_0, \bar t)\subset Q_R^\theta(x_0, t_0)$ we have that
 \begin{equation}
\left|\left\{(x,t)\in Q^\eta_R (x_0, \bar t): u(x,t)\ge \mu_+-\frac \omega{2^{s_1}}\right\}\right|
\le \left|\left\{ (x,t)\in Q^\eta_R(x_0, \bar t):u(x,t)\ge \mu_+-\frac \omega 2\right\}\right|\le (1-\nu)\left|Q^\eta_R(x_0, \bar t)\right|.
\label{e6.7}
 \end{equation}
 Further on we assume that \eqref{e6.1} holds.
\begin{lemma}\label{lem6.1}
Fix a cylinder $Q_R^\eta(x_0, \bar t)$. Suppose that \eqref{e6.7} holds.
There exists $t_\ast \in (\bar t-\eta R^p, \bar t-\frac {\nu \eta}2 R^p)$ such that
\begin{equation}
\left|\left\{x \in B_R(x_0):u(x,t_\ast)\ge \mu_+-\frac \omega{2^{s_1}}\right\}\right|\le \frac{1-\nu}{1-\frac\nu 2}\left|B_R(x_0)\right|.
\label{e6.8}
\end{equation}
\end{lemma}
\proof
Suppose not.  Then for all $t\in (\bar t-\eta R^p, \bar t-\frac{\nu \eta}2 R^p)$ there holds
\[\left|\left\{x\in B_R(x_0): u(x,t)\ge \mu_+-\frac\omega{2^{s_1}}\right\}\right|>\frac {1-\nu}{1-\frac{\nu}{2}}B_r(x_0).
\]
Hence
\begin{eqnarray*}&&\left|\left\{ (x,t)\in Q^\eta_R(x_0, \bar t):u(x,t)\ge \mu_+-\frac\omega{2^{s_1}}\right\}\right|\\
&\ge& \int\limits_{\bar t-\eta R^p}^{\bar t-\frac{\nu \eta}2R^p}\left|\left\{x\in B_R(x_0): u(x,t)\ge \mu_+-\frac\omega{2^{s_1}}\right\}\right| dt> (1-\nu)|Q_R^\eta(x_0, \bar t)|,
\end{eqnarray*}
which contradicts \eqref{e6.7}. \qed
\begin{lemma}\label{lem6.7}
There exists a number $s_1<s_2<s^\ast$, which depends only on known data, such that
\begin{equation}
\left|\left\{x\in B_R(x_0):u(x,t)\ge \mu_+-\frac \omega{2^{s_2}}\right\}\right|\le \left(1-\left(\frac \nu 2\right)^2\right)|B_R(x_0)|,
\label{e6.9}
\end{equation}
for all $t\in (\bar t-\frac \nu 2 \eta R^p, \bar t)$.
\end{lemma}
\proof
We use
Lemma~\ref{lem2.3} in the cylinder $B_R(x_0)\times (t_\ast, \bar t)$ with $k=\mu_+-\frac\omega{2^{s_1}}$, $\Psi$ defined by
\[
\Psi_+(u)=\ln_+ \frac{H^+_k}{H_k^+-(u-\mu_++\frac{\omega}{2^{s_2}})_++\frac{\omega}{2^{s_2}}},
\quad H^+_k=\esssup_{Q^\eta_R(x_0,\bar t)}(u-\mu_++\frac\omega{2^{s_1}})_+\le \frac\omega{2^{s_1}},
\]
and $\xi$ satisfying
$\ind_{B_{R(1-\sigma)}(x_0)}\le \xi \le \ind_{B_R(x_0)},\ |\n\xi|\le \frac2{\sigma R}$.
We can assume without loss that $H_k^-> \frac{\o}{2^{s_1+1}}$, since otherwise the assertion follows.
From Lemma~\ref{lem2.3} it follows that

\begin{eqnarray}
\nonumber
&& \int\limits_{B_{(1-\sigma)R}(x_0)}\Psi^2_+(u(x,t)) dx\le \int_{B_R(x_0)}\Psi^2_+ (u(x,t_\ast)) dx
+\gamma(\sigma R)^{-p} \iint\limits_{Q^\eta_R(_0, \bar t)}\Psi_+|\Psi_+^\prime(u)|^{2-p} dx dt\\
&+& \gamma\iint\limits_{Q^\eta_R(x_0, \bar t)}\Psi_+|\Psi_+^\prime(u)|^2 F_1(x) dx dt+\gamma\iint\limits_{Q^\eta_R(_0, \bar t)}\Psi_+|\Psi^\prime _+ (u)| F_2(x) dx dt.
\label{e6.10}
\end{eqnarray}
From the definition of $\Psi_+$ it follows that
\[
\Psi_+\le s_2\ln 2, \quad |\Psi^\prime(u)| \le \frac{2^{s_2}}\omega, \quad |\Psi_+^\prime(u)|^{2-p}\le \gamma\left(\frac\omega{2^{s_1}}\right)^{p-2}.\]

On the set $\left\{x\in B_R(x_0):u(x,t)\ge \mu_+-\displaystyle\frac \omega{2^{s_2}}\right\}$
we also have
$\Psi_+\ge (s_1-1) \ln2$.\\
Using \eqref{e6.8} we infer from \eqref{e6.10} that for all $t\in (t_\ast, \bar t)$
\begin{eqnarray}\nonumber
&&\left|\left\{x\in B_R(x_0):u(x,t)\le \mu_+-\frac \omega{2^{s_2}}\right\}\right|\\
\nonumber
&\le& (s_2-1)^{-2}\ln^{-2}2\int_{B_{(1-\sigma)R}(x_0)}\Psi^2_+(u(x,t)) dx +n\sigma|B_R(x_0)|\\
\nonumber
&\le& \left(\frac{s_2}{s_2-1}\right)^2 \frac{1-\nu}{1-\frac\nu 2}|B_R(x_0)|+n\sigma|B_R(x_0)|\\
&+&\gamma\frac{s_2}{(s_2-1)^2}\left\{\sigma^{-p}+\left(\frac{2^{s_2} \calF_1(2R)}\omega\right)^p+\left(\frac{2^{s_2}\calF_2(2R)}\omega\right)^{p-1}\right\}|B_R(x_0)|.
\label{e6.11}
\end{eqnarray}
First choosing $\sigma$ such that $n\sigma\le \frac 38=\nu^2$ and then $s_2$ such that
\[
\left(\frac{s_2}{s_2-1}\right)^2\le (1-\frac \nu 2)(1+\nu), \quad \gamma \frac{s_2}{(s_2-1)^2}(1+\sigma^{-p})\le \frac 38 \nu^2, \]
due to \eqref{e6.1} we obtain the required \eqref{e6.9} from \eqref{e6.11}.\qed

\medskip
Since inequality \eqref{e6.9} holds true for all cylinders $Q_R^\eta(x_0, \bar t)$, Lemma~\ref{lem6.1} implies the following assertion.

\begin{remark}\label{rem6.1}
For all $t\in (t_0-\frac12\theta R^p, t_0)$ the inequality
\begin{equation}\label{e6.12}
\left|\left\{ x\in B_R(x_0):u(x,t)\ge \mu_+-\frac \omega{2^{s_2}}\right\}\right|\le \left(1-(\frac \nu 2)^2\right)|B_R(x_0)|
\end{equation}
holds.
\end{remark}
\begin{lemma}\label{lem6.3}
For any $\nu\in (0,1)$ there  exists a number  $s^*$,  $s_2<s^\ast<\log_2b_2$, depending on the data only, such that
\begin{equation}\label{e6.13}
\left|\left\{(x,t)\in Q_R^\theta(x_0, t_0):u(x,t)\ge \mu_+-\frac\omega{2^{s^\ast}}\right\}\right|\le \nu\left|Q^\theta_R(x_0, t_0)\right|.
\end{equation}
\end{lemma}
The proof of Lemma \ref{lem6.3} is completely analogous to that of  Lemma \ref{lem5.1}.

Using Theorem \ref{thm1.4}
 with $\xi=\frac 1{2^{s^\ast}}, a=\frac12, \theta=\left(\frac {2^{s^\ast}}\omega\right)^{p-2}$ and $\nu$ defined by \eqref{e4.44}, from \eqref{e6.13} we obtain that
 \begin{equation}\label{e6.14}
 u(x,t)\le \mu_+-\frac\omega{2^{s^\ast+1}}=\mu_+-\sigma_1\omega\ \ \mbox {for a.a.}\   (x,t)\in Q^\theta_{\frac R2}(x_0, t_0).
 \end{equation}

Thus we have proved the following
 \begin{proposition}\label{prop6.2}
 Let the second alternative hold.
 Then either $\o\le b_2(R+\calF_1(2R)+\calF_2(2R))$, or
 \begin{equation}\label{e6.15}
 \essosc\limits_{Q_{\frac{R}2}^\theta(x_0, t_0)} u(x,t)\le (1-\sigma_1)\omega.
 \end{equation}
 \end{proposition}
 From Propositions~\ref{prop6.1},~\ref{prop6.2} in the same way as in \cite[Chapter III, Proposition 3.1]{DiB}
with the help of \cite[Lemma 8.23]{GT} we obtain:

\begin{proposition}\label{prop6.3}
For any $\varepsilon\in(0,1)$ and for all $\rho\le R$, there exist $\beta,\gamma>0$ and $\alpha\in(0,1)$, depending only on the data, such that
\begin{equation}
\essosc\limits_{Q(\rho, M)} u (x,t)\le\gamma\left(\frac \rho R\right)^\alpha \omega(R)+\gamma \calF_1(2\rho^\varepsilon R^{1-\varepsilon})+\gamma \calF_2(2\rho^\varepsilon R^{1-\varepsilon}),
\label{e6.16}
\end{equation}
where $Q(\rho, M)=B_\rho(x_0)\times(t_0-\beta M^{2-p}\rho^p, t_0), \ M=\esssup\limits_{\Omega_T}|u(x,t)|$.
\end{proposition}
This completes the proof of Theorem \ref{thm1.1}.

\section{Harnack inequality. Sketch of Proof of Theorem~\ref{Harnack}}
\label{harnack}

After we have proved Theorems~\ref{thm1.4}~and~\ref{thm1.5} the rest of the arguments
do not differ from  \cite{DiGV}. We give a short sketch here.

Let us consider the cylinder $Q_\tau=B_{\tau\rho}(x_0)\times \left(t_0-\frac{\tau^p\rho^p}{u_0^{p-2}},t_0\right)$, $u_0:=u(x_0,t_0)$.
%
Following Krylov-Safonov \cite{KS}
consider the equation
\[
\max_{Q_\tau}u(x,t)=u_0(1-\tau)^{-\b}
\]
where $\b>1$ is to be determined only depending on the data.
Let $\tau_0$ be the maximal root of the above equation and $u(\bar x ,\bar t)=u_0(1-\tau)^{-\b}$.
Let $\tilde Q=B_{\frac{1-\tau_0}2\rho}(\bar x)\times \left(\bar t-\left(\frac{1-\tau_0}2\right)^p\frac{\rho^p}{u_0^{p-2}},\,\bar t\right)$.
Since $\tilde Q\subset Q_{\frac{1+\tau_0}2}\subset Q_1$,
we have that
\[
\max_{\tilde Q} u\le \max_{Q_{\frac{1+\tau_0}2}} u\le 2^{\b}(1-\tau_0)^{-\b}u_0.
\]

{\it Claim~1.} There exists a positive number $\nu(\b)$ such that
\[
\left|\{(x,t)\in \tilde Q\,:\,u(x,t) \ge \frac12 (1-\tau_0)^{-\b}u_0 \} \right|>\nu(\b)|\tilde Q|.
\]
Indeed, in the opposite case we apply
Theorem~\ref{thm1.4} with the choices

$\mu_+=2^\b(1-\tau_0)^{-\b}u_0$, $\xi\o=(2^\b-\frac12)(1-\tau_0)^{-\b}u_0$, $a=\frac{2^\b-\frac34}{2^\b-\frac12}$.
The condition $u_0\ge B(\rho+\calF_1(2\rho)+\calF_2(2\rho))$ obviously implies that $\xi\o\ge B(\frac{1-\tau_0}2\rho+\calF_1({(1-\tau_0)}\rho)+\calF_2({(1-\tau_0)}\rho))$.
Therefore
we can conclude that
$u(\bar x,\bar t) \le \frac34(1-\tau_0)^{-\b}$ reaching
a contradiction which proves the claim.

{\it Claim~2.} (Analogue of \cite[Proposition~8.3]{DiGV}) For every $\nu_0\in(0,1)$ there exists a point $(y,s)\in\tilde Q$ and $\eta_0\in(0,1)$ and a cylinder
$Q_*=(y,s)+Q^\theta_{\eta_0 (1-\tau_0)\rho}\subset \tilde Q$ such that the inequality $u_0\ge B(\rho+\calF_1(2\rho)+\calF_2(2\rho))$ implies that
\[
|\{(x,t)\in Q_*\,:\,u(x,t)< \frac14(1-\tau_0)^{-\b}u_0|\le \nu_0|Q_*|.
\]
The proof is the same as in \cite{DiGV}. One writes down the energy inequality \eqref{e2.4} with $k=\frac12(1-\tau_0)^{-\b}u_0$ over coaxial cylinders $2\tilde Q$ and $\tilde Q$
and obtains the inequality
\[
\iint_{\tilde Q} |\n (u-k)_-|^pdxdt \le \g \frac{k^p}{R^p}|\tilde Q|,\quad R=\frac{1-\tau_0}2\rho.
\]
One only needs to note the estimates of the additional terms in the energy inequality \eqref{e2.4}. In the following
\eqref{e1.8b},\eqref{e1.9b} are used.
\[
\iint_{\tilde Q\cap\{u>k\}}f_1^\frac{p}{p-1} u^p \xi^p dxdt\le \g  \left(\frac{1-\tau_0}2\right)^p \rho^p u_0^{2-p} \left(\frac{1-\tau_0}2\rho\right)^{n-p}\calF_1(2\rho)^p
\le \g \frac{k^p}{R^p}|\tilde Q|,
\]
\[
\iint_{\tilde Q}(u-k)_-f_2\xi^p dxdt\le \g k \left(\frac{1-\tau_0}2\right)^p \rho^p u_0^{2-p} \left(\frac{1-\tau_0}2\rho\right)^{n-p}\calF_2(2\rho)^{p-1}
\le \g \frac{k^p}{R^p}|\tilde Q|.
\]

The rest of the proof of Claim~2 is the same as in  \cite[Proposition 8.3]{DiGV} and is based on Lemma~\ref{lem2.1}, which in turn relies on \cite{DiGV-1}.

As in \cite{DiGV}, by Theorem~\ref{thm1.4}
we obtain that there exist $(y,s)\in\tilde Q$
and $\eta_0\in (0,1)$ such that
\[
u(x,s)\ge \frac{u_0}{8}(1-\tau_0)^{-\b}\quad \text{for}\ |x-y|\le r:=\eta_0\frac{1-\tau_0}2\rho.
\]
Then an application of Theorem~\ref{thm1.5} yields
that if
\[
\frac{u_0}8(1-\tau_0)^{-\b} \ge B(r+\calF_1(2r)+\calF_2(2r))
\]
then
\[
u(x,t)\ge \sigma \frac{u_0}8(1-\tau_0)^{-\b}\quad \text{for}\ |x-y|\le 2r, \ s+\left(\frac{u_0}8(1-\tau_0)^{-\b}\right)^{2-p}b_1 r^p \le t\le
s+\left(\frac{u_0}8(1-\tau_0)^{-\b}\right)^{2-p}b_2 r^p.
\]
After iteration for $j=1,2,3,\dots$ we have either
\[
\sigma^{j-1} \frac{u_0}8(1-\tau_0)^{-\b} \le B(\sigma^j r+\calF_1(2\sigma^j r)+\calF_2(2\sigma^j r))
\]
or
\[
u(x,t)\ge \sigma^j \frac{u_0}8(1-\tau_0)^{-\b}\quad \text{for}\ |x-y|\le 2^j r,
\]
\[
\ t_j^{(1)}=t_{j-1}^{(1)}+\left(\frac{u_0}8(1-\tau_0)^{-\b}\sigma^j\right)^{2-p}b_1 r^p \le t\le
t_{j-1}^{(2)}+\left(\frac{u_0}8(1-\tau_0)^{-\b}\sigma^j\right)^{2-p}b_2 r^p=t_j^{(2)}.
\]
Choosing $j$ such that $2^j \eta_0 \frac{1-\tau_0}2=2$ and $\b$ by the condition $2^\b\sigma=1$ we complete the proof.
 (see \cite[Section~8]{DiGV} for details).

\begin{small}
\section*{Acknowledgments}
The authors acknowledge support of the Royal Society through the
International Joint Project Grant 2006/R1.
\end{small}
\begin{small}

\end{small}


\end{document}